\documentclass[11pt,reqno]{amsart}

\usepackage{euscript}
\usepackage{epsfig}
\usepackage{amssymb}

\theoremstyle{plain}
\newtheorem{theorem}{Theorem}[subsection]

\newtheorem{proposition}[theorem]{Proposition}

\theoremstyle{definition}

\newcommand{\no}{\noindent }
\newcommand{\bc}{{\mathbb C}}
\newcommand{\bp}{{\mathbb P}}
\newcommand{\bz}{{\mathbb Z}}
\newcommand{\bn}{{\mathbb N}}
\newcommand{\bq}{{\mathbb Q}}
\newcommand{\vp}{\varphi}
\newcommand{\vpt}{\varphi_{2\ell+3}}
\newcommand{\vpf}{\varphi_{2\ell+5}}
\newcommand{\vpo}{\varphi_{2\ell+1}}
\newcommand{\vpm}{\varphi_{2\ell-1}}
\newcommand{\vps}{\varphi_{2\ell+7}}
\newcommand{\vptw}{\varphi_{2\ell+2}}
\newcommand{\blambda}{\boldsymbol{\lambda}}

\newcommand{\calA}{{\mathcal{A}}}
\newcommand{\calB}{{\mathcal{B}}}
\newcommand{\calG}{{\mathcal{G}}}

\newcommand{\comp}{{\circ}}

\begin{document}

\title{On rational cuspidal projective plane curves}
\author{J. Fern\'andez de Bobadilla,
I. Luengo-Velasco, A. Melle-Hern\'andez, and A. N\'emethi}
\address{Department of Mathematics, University of Utrecht, Postbus 80010, 3508TA Utrecht, The Netherlands}
\address{Facultad de Matem\'aticas\\ Universidad Complutense\\ Plaza de Ciencias\\ E-28040, Madrid, Spain}
\address{Department of Mathematics\\Ohio State University\\Columbus,
OH 43210,USA; and R\'enyi Institute of Mathematics, Budapest,
Hungary.}
\address{}
\email{bobadilla@math.uu.nl} \email{iluengo@mat.ucm.es}
\email{amelle@mat.ucm.es} \email{nemethi@math.ohio-state.edu;
nemethi@renyi.hu}
\thanks{The first author is supported by the Netherlands Organization for Scientific Research (NWO). The first three authors are partially
supported by BFM2001-1488-C02-01. The forth
 author is partially supported by NSF grant DMS-0304759.}




\maketitle
\pagestyle{myheadings}
\markboth{{\normalsize J. Fern\'andez de Bobadilla, I. Luengo, A. Melle-Hern\'andez, A.
N\'emethi}}{
{\normalsize Rational cuspidal plane curves}}

{\small

\section{Introduction}

Let $C$ be an irreducible   projective plane curve
in the complex projective space $\bp^2$   with singular points $\{p_i\}
_{i=1}^\nu$. It is a very interesting, and still open problem, to characterize
the local
embedded topological types of local singular germs $(C,p_i)\subset (\bp^2,p_i)
$ which can be realized by such
a projective curve $C$ of degree $d$.
This remarkable
 problem is not only important for its own sake, but it is also
connected with crucial properties, problems and conjectures in the theory of
open surfaces.

For instance, the open surface $\bp^2\setminus C$ is
$\bq$-acyclic if and only if $C$ is a rational cuspidal curve. On
the other hand, the \emph{rigidity conjecture} proposed by
 Flenner and Zaidenberg  in \cite{FZ0}
(and supported by many examples, see \cite{FZ0}, \cite{FZ1},
\cite{FZ2}, \cite{Fenske1}) says that every $\bq$-acyclic affine
surfaces $Y$ with logarithmic Kodaira dimension
$\bar{\kappa}(Y)=2$ must be rigid. (E.g., if $C$ has at least
three cusps then $\bar{\kappa}(\bp^2\setminus C)=2$, cf.
\cite{Wak}.)  This conjecture for $Y=\bp^2\setminus C$ would
imply the projective rigidity of  the curve $C$ in the sense that
every equisingular deformation of $C$ in $\bp^2$ would be
projectively equivalent to $C$.

Among other interesting related open problems we mention: every
rational cuspidal curve can be transformed by a Cremona
transformation into a line (this is the Coolidge-Nagata problem,
see \cite{Coo},\cite{Na}); or, the determination of the maximal
number of cusps among all the rational cuspidal plane curves
(proposed by F. Sakai in \cite{open}) -- this number  is expected
to be small (the maximal known by the authors is four). In a
recent preprint by K.~Tono \cite{Tono2} it is proved that the
maximal number is strictly less than nine. Finally, we also list
the conjectured numerical inequality (2.3) \cite{Or} of
S.Yu.~Orevkov. [In forthcoming articles we plan
to run the machinery of the
present manuscript for these open problems as well -- at least
when $\bar{\kappa}(\bp^2\setminus C)<2$.]

The first-mentioned  `characterization  problem' (on the
realization of prescribed topological types of singularities) has
a long and rich history providing many interesting compatibility
properties connecting local invariants of the germs
$\{(C,p_i)\}_i$ with some global invariants of $C$ (like its
degree, or the log-Kodaira dimension of $\bp^2\setminus C$,
etc.). Most of the compatibility properties are mere identities
or inequalities (see e.g. the Matsuoka-Sakai inequality, or
Orevkov's sharp inequality, both being consequences of the
logarithmic version of the Bogomolov-Miyaoka-Yau inequality). The
more complex Varchenko's compatibility property is provided by the
semicontinuity of the spectral numbers of isolated hypersurface
singularities,  and it consists of a finite set of inequalities.

Our goal is to propose a new compatibility
property -- valid for rational cuspidal curves $C$ --
which seems to be rather  powerful
(see section 3 for some comparison with  some classical criterions).
Its formulation is surprisingly very elementary.
Consider a collection $(C,p_i)_{i=1}^\nu$ of locally irreducible
plane curve singularities (i.e. cusps), let $\Delta_i(t)$ be the
characteristic polynomials of the monodromy action associated with
$(C,p_i)$, and $\Delta(t):=\prod_{i}\Delta_i(t)$.
Its degree is $2\delta$, where $\delta$ is the sum of the delta-invariants
of the singular points.
Then $\Delta(t)$ can be written as
$1+(t-1)\delta+(t-1)^2Q(t)$
for some polynomial $Q(t)$. Let $c_l$ be the coefficient of $t^{(d-3-l)d}$
in $Q(t)$ for any $l=0,\ldots, d-3$.

\vspace{2mm}

\noindent {\bf Conjecture.} {\em Let  $(C,p_i)_{i=1}^\nu$ be a
 collection of local plane curve  singularities, all of them locally
irreducible,  such that
$2\delta=(d-1)(d-2)$ for some integer $d$. Then if $(C,p_i)_{i=1}^\nu$
can be realized as the local singularities of a degree $d$
(automatically rational and cuspidal)
projective plane curve of degree $d$ then}
\begin{equation*}c_l\leq (l+1)(l+2)/2 \ \ \mbox{for all \ $l=0,\ldots, d-3$.}
\tag{$*_l$}\end{equation*}
In fact, the integers $N_l:=c_l-(l+1)(l+2)/2$ \ are symmetric: $N_l=N_{d-3-l}$;
and  $N_0=N_{d-3}=0$ automatically.
Moreover, there is a surprising phenomenon in the above conjecture:

\vspace{2mm}

{\em If $\nu=1$, then the conjecture is true if and only if in all the
inequalities ($*_l$), in fact, one has equality (cf. \ref{nuone}). }

\vspace{2mm}

The conjecture can be reformulated in the language of the semigroups of
the germs $(C,p_i)$ (and the degree $d$) as well. E.g., if $\nu=1$,
then collection of vanishings of all the coefficients $N_l$ can be
described by a very precise and mysterious distribution of the
elements of the semigroup of the unique singular point with respect to the
intervals $I_l:=(\,(l-1)d,ld\,]$ (cf. section 3):

\vspace{2mm}

{\em If $\nu=1$, the conjecture predicts that
there are exactly $\min\{l+1,d\}$ elements of the semigroup in $I_l$
for any $l>0$. }

\vspace{2mm}

The idea (and the main motivation) of the above conjecture  came from the
Seiberg-Witten invariant  conjecture formulated for normal surface
singularities by Nicolaescu and the forth author \cite{[51]}, respectively
from the set of counterexamples for this conjecture provided by superisolated
singularities \cite{SI} found by the last three authors.
For the completeness of the presentation,
this is described in short in section 2, but we emphasize that the present
article is independent of the techniques of \cite{[51]}, in particular
involves no Seiberg-Witten theory. (On the other hand, we are witnesses
of a mysterious connection whose deeper understanding would be a wonderful
mathematical goal.)

As supporting evidence, in
the body of the paper we prove the following result:

\vspace{2mm}

\noindent
{\bf Theorem 1.} {\em If the logarithmic Kodaira dimension
$\bar{\kappa}:=\bar{\kappa}(\bp^2\setminus C)$ is $\leq 1$,
then the above conjecture
is true. In fact, in all these cases $N_l=0$ for any $l=0,\ldots, d-3$
(regardless of $\nu$). }

\vspace{2mm}

The proof of Theorem 1 consists of  many steps. Its structure is the following.

\vspace{2mm}

(a) If $\bar{\kappa}=-\infty$ then $\nu=1$ by \cite{Wak}. Moreover, all these
curves are classified by H. Kashiwara \cite{Kash}. The family contains as an
important subfamily the Abhyankar-Moh-Suzuki ($AMS$) curves.

In this case, our proof runs as follows: first we verify the
vanishing of the coefficients $N_l$ for the $AMS$ curves (section 4) --
this corresponds to the ${\mathcal F}_I$ family in Kashiwara's classification.
The case of the other subfamily (the ${\mathcal F}_{II}$ curves)
is treated in section 6. The proof is based on an induction where the starting
point is given by some unicuspidal curves which have only one
Puiseux pair. This case is completely resolved in section 5.

(b) The case $\bar{\kappa}=0$ cannot occur by a result of
Sh.~Tsunoda \cite{Tsu1},
see also the paper by Orevkov \cite{Or}.

(c) If $\bar{\kappa}=1$ then by a result of Wakabayashi \cite{Wak}
one has $\nu\leq 2$. In the case  $\nu=1$, K. Tono writes the
possible equations of the curves \cite{KeitaTono}. (Notice that
Tsunoda' classification  in \cite{Tsu} is incomplete.) We verify
the conjecture for them in section 7. On the other hand, by
another  result of Tono \cite{To}, the case $\nu=2$ corresponds
exactly to the Lin-Zaidenberg bicuspidal rational plane curves.
For them we verify the conjecture in section 8.

\vspace{2mm}

We wish to emphasize that in the process of the verification of
the conjecture, in fact, {\em we list all the possible local
topological types} (e.g. Eisenbud-Neumann splice diagrams) {\em
of local plane curve singularities, together with the degrees
$d$, which can be geometrically realized with $\bar{\kappa}\leq
1$}.  But, from the point of view of the proof of Theorem 1, this
information is far to be enough to verify the conjecture. Even if
one knows the local topological types (e.g. the local resolution
graphs, or even better, the generators of the corresponding
semigroups), one needs sometimes additional rather involved
arithmetical arguments to complete the proof.

It is important to notice that in the $\bar{\kappa}=2$ case there are
(infinitely many, and for arbitrarily large $d$) examples when {\em
not all}  the
coefficients $N_l$ vanish. The complete picture for curves
with $d\leq 6$ is provided in section 2,
together with some  additional introductory examples  with larger $d$.
The following list provides those cases
which are verified in the present article.

\vspace{2mm}

\noindent {\bf Theorem 2.} {\em If $\bar{\kappa}=2$, then in the following
cases the conjecture is true:

(a) $d\leq 6$;

(b) $C$ is unicuspidal with one Puiseux pair (see section 5);

(c) Orevkov's curves $C_{4k}$ and $C_{4k}^*$ \cite{Or} (see section 9).

(d) $d$ even, $\nu=2$, and  multiplicity sequences
 $[d-2], [2_{d-2}]$ (see \ref{dbig}).}

\vspace{2mm}

\noindent See also  remark \ref{3.12}, valid for $\nu=1$.

\vspace{2mm}

We thank J. Stevens for drawing our attention to an error in the original
version of \ref{varch}, and K. Tono for providing us with extremely helpful
information.

\section{The Main Conjecture. Motivation and first comments.}

For the convenience of the reader, we recall some notations and classical
properties of plane curve singularities, which will be intensively used in the
sequel.

\subsection{Invariants of germs of irreducible plane curve
 singularities.}\label{a1}
To encode the topology of  a germ of an irreducible plane   curve singularity
$(C,0)\subset (\bc^2,0)$ several sets of invariants can be used: Puiseux pairs, Newton pairs,
(minimal) embedded resolution graph,  Eisenbud-Neumann splice diagram, semigroup, etc.
We mainly use the  Eisenbud-Neumann splice diagram
(cf. \cite{EN}, page 49). If the germ $(C,0)$ has $g$  Newton pairs
$\{(p_k,q_k)\}_{k=1}^g$
with gcd$(p_k,q_k)=1$ and $p_k\geq 2$ and $q_k\geq 2$
(and by convention, $q_1>p_1$),
define the integers $\{a_k\}_{k=1}^g$ by
$a_1:=q_1$ and $a_{k+1}:=q_{k+1}+p_{k+1}p_ka_k$ for  $k\geq 1$.
Then its  Eisenbud-Neumann splice diagram decorated by
the numerical data $\{(p_k,a_k)\}_{k=1}^g$ has the following  shape:

\begin{picture}(400,65)(0,10)
\put(50,60){\circle*{4}}
\put(100,60){\circle*{4}}
\put(150,60){\circle*{4}}
\put(250,60){\circle*{4}}
\put(300,60){\circle*{4}}
\put(100,20){\circle*{4}}
\put(150,20){\circle*{4}}
\put(250,20){\circle*{4}}
\put(300,20){\circle*{4}}
\put(350,60){\makebox(0,0){$\tilde{C}$}}
\put(92,65){\makebox(0,0){$a_1$}}
\put(142,65){\makebox(0,0){$a_2$}}
\put(240,65){\makebox(0,0){$a_{g-1}$}}
\put(292,65){\makebox(0,0){$a_g$}}
\put(108,50){\makebox(0,0){$p_1$}}
\put(158,50){\makebox(0,0){$p_2$}}
\put(262,50){\makebox(0,0){$p_{g-1}$}}
\put(308,50){\makebox(0,0){$p_g$}}
\put(200,60){\makebox(0,0){$\cdots$}}
\put(50,60){\framebox(125,0){}}
\put(225,60){\framebox(75,0){}}
\put(100,20){\framebox(0,40){}}
\put(150,20){\framebox(0,40){}}
\put(250,20){\framebox(0,40){}}
\put(300,20){\framebox(0,40){}}
\put(300,60){\vector(1,0){30}}
\end{picture}

The characteristic polynomial $\Delta_{(C,0)}(t)$ of the
monodromy acting on the first homology of the Milnor fiber of the singularity
can be computed by A'Campo's formula \cite{AC} from the splice diagram.
If we define
\begin{equation*}
\begin{array}{ll}
{\beta_k}:=a_kp_kp_{k+1}\cdots p_g & \mbox{for $1\leq k\leq g$};\\
{\bar{\beta}_0}:=p_1p_2\cdots p_g; & \\
{\bar{\beta}_k}:=a_kp_{k+1}\cdots p_g & \mbox{for $1\leq k\leq g$}, \end{array}
\end{equation*}
then  $\Delta_{(C,0)}$
is given by:
\begin{equation*}
\Delta_{(C,0)}(t)= \frac{(t-1)\prod_{1\leq k\leq g} (t^{\beta_k}-1)}
{\prod_{0\leq k\leq g} \, (t^{\bar{\beta}_k}-1)}.
\end{equation*}
The polynomial $\Delta_{(C,0)}$ is a complete (embedded)
topological invariant of the germ $(C,0)\subset (\bc^2,0)$,
similarly as the semigroup $\Gamma_{(C,0)}\subset \bn$ generated
by all the possible intersection multiplicities $i(\{h=0\},C)$ at
$0$ for all $h\in {\mathcal O}_{(\bc^2,0)}$.  The degree
$2\delta_{(C,0)}$ of $\Delta_{(C,0)}$ is the conductor of the
singularity, where the delta-invariant $\delta_{(C,0)}$ is the
cardinality of the finite set $\bn\setminus \Gamma_{(C,0)}$.

By \cite{CDG}, $\Delta_{(C,0)}(t)=(1-t)\cdot L(t)$, where $L(t)
:=\sum_{k\in \Gamma_{(C,0)}} t^k$ is the
Poincar\'e series of $\Gamma_{(C,0)}$.
 In fact, the minimal set of generators of
$\Gamma_{(C,0)}$ consists of the
$g+1$ elements $\bar{\beta}_i $  ($0\leq i\leq g$).
It is also known that each element $\gamma\in \Gamma_{(C,0)}$
can be represented in a unique way in
the form $\gamma=k_0\bar{\beta}_0+\sum_{1\leq j\leq g} k_j\bar{\beta}_j$
with $k_0\geq 0$ and $0\leq k_j \leq p_j-1$ for $1\leq j\leq g$, see \cite{TZ}.

\subsection{Motivation.
The Seiberg-Witten invariant conjecture ($SWC$).}\label{2.2} In
\cite{[51]}
L. Nicolaescu and the forth  author formulated the following conjecture
(as a generalization of the ``Casson invariant conjecture'' of
Neumann and Wahl \cite{NW}):

\vspace{2mm}

 {\em If the link of a normal surface singularity $(X,0)$ is a rational
homology sphere then the geometric genus $p_g$ of $(X,0)$ has an
``optimal'' topological upper bound. Namely,
\begin{equation*}p_g \leq  {\bf sw}(M)-(K^2+s)/8.\tag{$SWC$}\end{equation*}
Moreover, if $(X,0)$ is a  ${\mathbb Q}$-Gorenstein (e.g.
hypersurface) singularity then in ($SWC$) the equality holds.}

\vspace{2mm}

Here, ${\bf sw}(M)$ is the (topological, or `modified') Seiberg-Witten
invariant of the link $M$ of $(X,0)$ associated with its canonical $spin^c$
structure, $K$ is the canonical cycle associated with a fixed resolution
graph $G$ of $(X,0)$, and $s$ is the number of vertices of $G$
(see \cite{[51]} for more details).

\vspace{2mm}

The ($SWC$)-conjecture  was verified successfully for many different families,
see e.g. \cite{[51],[52],[55],NOSZ}. But the last three authors in \cite{SI}
found some counterexamples based on superisolated singularities.
This class ``contains'' in a canonical way the theory of complex projective
plane curves, a fact which is crucial in the next discussion.

Hypersurface superisolated singularities
were introduced in \cite{Ignacio} by the second author and achieved
the reputation of being a distinguished class of singularities and
source of interesting examples and counterexamples.
A hypersurface singularity $f:(\bc^3,0)\to (\bc,0)$, $f=f_d+f_{d+1}+\cdots$
(where $f_j$ is homogeneous of degree $j$) is superisolated if the
projective plane curve $C:=\{f_d=0\}\subset \bp^2$ is reduced with
isolated singularities $\{p_i\}_{i=1}^\nu$, and these points are not situated
on the projective curve $\{f_{d+1}=0\}$. In this case the embedded
topological type (and the equisingular type) of $f$  depends only
on the curve $C$. Notice also that the link of $f$ is a
rational homology sphere if the curve $C$ is rational and cuspidal
(i.e. if all the germs $(C,p_i)$ are locally irreducible).
{\em In the sequel we also will assume these two facts. }

In \cite{SI} the authors have shown
that some hypersurface superisolated singularities
with $\nu=\#Sing(C)\geq 2$ do not satisfy the above Seiberg-Witten
invariant conjecture. Moreover,  in all the counterexamples
   $p_g>{\bf sw}(M)-(K^2+s)/8$ (contrary to the
inequality predicted by the general conjecture \ref{2.2} !).
 On the other hand, even after an intense search
of the existing cases, the authors were not able to find any
counterexample with $\nu=1$. (In fact, one of the goal of the
present article is to explain, at least partially, these phenomenons, cf.
\ref{nuone}.)

In the next paragraphs, we plan to reformulate ($SWC$)
relating with an even deeper conjectural property.
We denote by $\Delta_i$ the characteristic polynomial of
$(C,p_i)\subset (\bp^2,p_i)$, set  $\Delta(t) :=\prod_i\Delta_i(t)$ and
$2\delta:=\deg\Delta(t)$.
By the rationality of $C$ one has   $(d-1)(d-2)=2\delta$.
Clearly, $\delta $ is the sum of the delta-invariants
of the germs $(C,p_i)_{i=1}^\nu$.
Notice also that (see e.g. \cite{[52]} (4.3)):
\begin{equation*}
\Delta(1)=1 \ \ \mbox{and} \ \  \Delta'(1)=\delta.
\tag{1}\end{equation*}
In the next discussion it is convenient to introduce the
polynomials $P$ and $Q$ by
$$\Delta(t)=1+(t-1)P(t)=1+(t-1)\delta+(t-1)^2Q(t).$$
The corresponding coefficients are denoted by
\begin{equation*}
P(t)=\sum_{l=0}^{2\delta-1}a_kt^k, \ \ \mbox{and}\ \
Q(t)=\sum_{l\nmid d}b_lt^l+\sum_{l=0}^{d-3}c_lt^{(d-3-l)d}.\tag{2}
\end{equation*}
The following definition helps to describe the properties of these
coefficients.
\subsection{Definition/Lemma.}\label{deflem} {\em We say that a polynomial
$D(t)=\sum_l\alpha_lt^l$ has the ``negative distribution property'' if
$\sum_{0\leq l\leq k}\alpha _l\leq 0$ for any $k\geq 0$.
Notice that if $D(t)=N(t)(1-t)$ for some other polynomial $N(t)$,
then $D(t)$ has the negative distribution property if and only if
all the coefficients of $N(t)$ are nonpositive. }

\vspace{2mm}

\noindent For the superisolated singularity  $(X,0):=\{f=0\}$
those invariants which appear in the formula ($SWC$) (cf.
\ref{2.2}) are computed  in \cite{SI} as follows (here ${\bf
sw}(M)$ is computed by the Reidemeister-Turaev torsion normalized
by the Casson-Walker invariant) :
\begin{equation*}
\left\{ \begin{array}{l}
K^2+s=-(d-1)(d^2-3d+1); \ \ \
p_g=d(d-1)(d-2)/6;\ \ \mbox{and} \\ \\
{\displaystyle
{\bf sw}(M)=\frac{1}{d}\sum_{\xi^d=1\not=\xi}\ \frac{\Delta(\xi)}{(\xi-1)^2}
 +\frac{1}{2d} \Delta(t)''(1) -\frac{\delta(6\delta-5)}{12d}.
}
\end{array}\right.
\tag{3}\end{equation*}

\vspace{2mm}

\noindent
Consider now the (a priori) rational function
\begin{equation*}R(t):= \frac{1}{d}\sum_{\xi^d=1}
\frac{\Delta(\xi t)}{(1-\xi t)^2}-\frac{1-t^{d^2}}{(1-t^d)^3}.\tag{4}
\end{equation*}

\subsection{Proposition.}\label{abcde} {\em With the above notations
one has:

(a) $R(t)$ can be written in the form $D(t^d)/(1-t^d)$, where
$$D(t)=(d-2)t+(d-3)t^2+\cdots+t^{d-2}-\sum_{k=0}^{2\delta-1}a_kt^{\lceil k/d
\rceil}$$
$$=\sum_{k\geq 0}(1-a_k)t^{\lceil k/d\rceil}-\frac{1-t^d}{(1-t)^2}.
\hspace{2cm}$$

(b) In (a), some combinations of the coefficients $\{a_k\}_k$
give no contribution in $R(t)$. Indeed, if one writes $P(t)=\delta+(t-1)Q(t)$
as above (cf. (2)), then the $\{b_l\}_l$ coefficients have no effect in $R$.
More precisely,
$R(t)=N(t^d)$, where
$$N(t)=\sum_{l=0}^{d-3} \Big( c_l-\frac{(l+1)(l+2)}{2}\Big)\, t^{d-3-l}.$$

(c) $N(t)$ is a symmetric polynomial (i.e.  $N(t)=t^{d-3}\cdot N(1/t)$) with
integral coefficients and with $N(0)=0$.

(d) $$R(1)={\bf sw}(M)-\frac{K^2+s}{8}-p_g.$$

\noindent In particular, the following facts also hold:

(e) $R(t)\equiv 0$ if and only if $D(t)\equiv 0$ if and only if
$N(t)\equiv 0$. In this case evidently ($SWC$) is true (with equality).

(f) $D(t)$ has the negative distribution property if and only if
all the coefficients of $R(t)$ (or, equivalently, of $N(t)$)
are nonpositive. In this case $p_g\geq {\bf sw}(M)-(K^2+s)/8$.
}

\begin{proof}  Write
$\Delta(t)/(1-t)^2=1/(1-t)^2-P(t)/(1-t).$
By an elementary computation
$$\frac{1}{d}\sum_{\xi^d=1}
\frac{1}{(1-\xi t)^2}=\frac{1+(d-1)t^d}{(1-t^d)^2}.$$
For the second contribution consider the monomials $a_kt^k$ of $P(t)$.
Express any fixed $k$ in the form $k=dn+r$ with $r=1,2,\ldots,d.$
Then
$$
\frac{1}{d}\sum_{\xi^d=1} \frac{\xi^k t^k}{1-\xi t}=\frac{t^{k}}{d}
\sum_{\xi^d=1} \frac{\xi^r}{1-\xi t}=
\frac{t^k}{d(1-t^d)}\sum _{\xi^d=1} (\, 1+\xi t+\cdots+
(\xi t)^{d-1})\xi^r =\frac{t^{d(n+1)}}{1-t^d}.$$
Notice also that
\begin{equation*}
\frac{1-t^d}{(1-t)^2}
=1+2t+\cdots+ (d-1)t^{d-2}+d(t^{d-1}+t^d+\ldots),\end{equation*}
and
$$
\frac{1+(d-1)t}{(1-t)}=1+d(t+t^2+\cdots)=\sum_{k\geq 0}\, t^{\lceil
 k/d\rceil},$$
which ends the proof of (a). The proof of (b) is similar. (c) follows from (b)
and from (4) using the fact that $\Delta$
 is monic and symmetric. For (d)  write
$$R(t)=
\frac{1}{d}\sum_{\xi^d=1\not=\xi}
\frac{\Delta(\xi t)}{(1-\xi t)^2}+
\frac{(1+t+\cdots +t^{d-1})^3\Delta(t)-d(1+t+\cdots +t^{d^2-1})}
{d(1-t)^2(1+t+\cdots +t^{d-1})^3}.$$
Then compute $\lim_{t\to 1}R(t)$ by the L'Hospital rule, then use (1) and (3).
\end{proof}

The starting point of the present article is the following conjecture.

\subsection{Conjecture.}\label{2.4}  {\em Assume
that $C$  is a rational cuspidal projective plane curve as above. Then the
property (f) of Proposition \ref{abcde} is always true. Namely, $N(t)$
(or equivalently, $R(t)$) has only nonpositive coefficients.
In other words, $D(t)$ has the negative distribution property. }

\vspace{2mm}

In fact, one can ask about the stronger property \ref{abcde}(e),
i.e. about the vanishing of $N(t)$ (or equivalently, of $R(t)$).
Obviously, this does not hold in general: all the examples of
\cite{SI} with $p_g> {\bf sw}(M)-(K^2+s)/8$ have $R(t)\not\equiv
0$ (cf. also with \ref{dsmall}). On the other hand, using
\ref{abcde}(d) one gets the following: {\em if one can verify
independently the inequality $p_g\leq  {\bf sw}(M)-(K^2+s)/8$
then the above Conjecture \ref{2.4} is true if and only if, in
fact, $R(t) \equiv 0$}. This is the case, e.g., if $\nu=1$:

\subsection{Theorem.}\label{nuone} {\em Assume that $C$ is an
unicuspidal  rational curve. Then Conjecture \ref{2.4} is true if and only if
$R(t)$ is identical zero.}

\vspace{2mm}

We will give two proofs of this theorem. The first one, given below,
is based on a  geometrical
result which reveals the connection of Conjecture~\ref{2.4}
with links of surface singularities and topological invariants of 3-manifolds.
 The second (see section   \ref{secgenus}) is
rather elementary and surprising.

\begin{proof}
Recall that a good resolution graph of a normal surface singularity
is called $AR$ (almost rational) if the graph has at least one  vertex
with the following property: if one replaces  the decoration
(self-intersection) of the corresponding vertex by a smaller
(more negative) integer then one gets a rational graph (in the sense
of Artin); for more details see \cite{NOSZ}. The point is that if a normal
surface singularity with rational homology sphere link admits an
$AR$ good resolution graph then it satisfies the inequality
$p_g\leq {\bf sw}(M)-(K^2+s)/8$, cf. \cite{NOSZ} (9.6).
[Here we have to mention that \cite{NOSZ} verifies the inequality for the
Ozsv\'ath-Szab\'o `Euler-characteristic' ${\bf sw}^{OSZ}(M)$,
but, by \cite{Rus},
this agrees with the Reidemeister-Turaev sign refined torsion
normalized by the Casson-Walker invariant, used here and in \cite{[51]}.]
Therefore,
it is enough to verify that $\{f=0\}$ admits an $AR$ good
resolution graph. Recall that the minimal good resolution graph $G$ of
$\{f=0\}$ can be obtained from the  minimal good embedded
resolution graph $\Gamma$ of the unique singular point $(C,p)$ of $C$
by gluing an extra vertex $w$ by a unique edge to the unique $-1$
vertex $v$ of $\Gamma$.
The decoration of $w$ (as any decoration)
 is negative, say $-r$, its precise value is not
essential here -- for details see e.g. \cite{SI}. Now, let $G_s$ be the
graph obtained from $G$ by changing
the decoration $-1$ of $v$ by $-2$. Then $G_s$ is the resolution graph
of a sandwiched
singularity. Indeed, consider again $\Gamma$.
Blow up a generic point on the curve corresponding to $v$,
and then blow up $r-1$ generic points on the new curve. In this way one
gets a non-minimal resolution graph of $(C,p)$
which contains the graph  $G_s$ as a subgraph.
\end{proof}

See also Example \ref{e3} for some additional comments about \ref{nuone}.

\subsection{Examples ($d$ small; $\nu\geq 1$).}\label{dsmall}
Below, any singular point $(C,p_i)$ will be identified by its multiplicity
sequence. Since the number of occurrences of the multiplicity $1$
in the multiplicity sequence equals the last multiplicity greater than $1$,
we omit the multiplicity $1$: we denote such a sequence by
$ [m_1,\ldots,m_l]$ where $ m_1 \ge m_2 \ge \cdots \ge m_l > m_{l+1} = 1$
for a suitable $l \ge 1$. In fact, we
will write $[\hat{m}_{1_{r_1}},\ldots,\hat{m}_{k_{r_k}}]$
for a multiplicity sequence which means that
the multiplicity $\hat{m}_i$ occurs $r_i$ times for
$i=1,\ldots,k$. For example, $[4_2,2_3]$ means $[4,4,2,2,2,1,1]$.

For the classification of the cuspidal rational curves with
degree $d\leq 5$, see e.g. the book of Namba \cite{Namba};
for the classification of multiplicity sequences of
rational cuspidal plane curves of degree $6$
see e.g.  Fenske's paper \cite{Fenske}.

If $d=3$, then $C$ has a unique singularity of type $[2]$.
If $d=4$, then there are four possibilities; the corresponding
multiplicity sequences of the singular points are
$[3]$; $[2_3]$; $[2_2],[2]$ and $[2],[2],[2]$.
By a verification (or by  \ref{abcde}(c)),
in all  these cases $N(t)\equiv 0$, hence
the conjecture is true.

For  $d=5$ and $d=6$ the next table shows all the possible multiplicity
sequences together with the polynomials $N(t)$. For all these cases
the conjecture is also true.

\vspace{2mm}

\begin{minipage}[t]{2in}
    \begin{tabular}{l|l|l|l}
           $d$ & $\nu$   & type of cusp       & $-N(t)$ \\ \hline
            5  &    1        & $[4]$           &   0         \\
            5 &     1       & $[2_6]$      &    0        \\
            5 &     2       & $[3,2]$ , $[2_2]$  & 0   \\
           5  &  2         & $[3]$ , $[2_3]$  &  $2t$  \\
          5   &   2        & $[2_2]$ , $[2_4]$ & $2t$ \\
       5         & 3      &  $[3]$ , $[2_2]$ , $[2]$ & $2t$   \\
        5 &       3      & $[2_2]$ , $[2_2]$ , $[2_2]$ & $6t$  \\
       5 & 4        & $[2_3]$ , $[2]$ , $[2]$, $[2]$ & $8t$   \\
     \end{tabular}
\end{minipage} \hspace{1cm}
\begin{minipage}[t]{2in}
    \begin{tabular}{l|l|l|l}
           $d$ & $\nu$   & type of cusp       & $-N(t)$ \\ \hline
     6 &    1    & $[5]$           &  0 \\
      6 & 1   & $[4,2_4]$      & 0 \\
      6 & 1    & $[3_3,2]$      &   0 \\
      6 & 2   & $[3_3]$  , $[2]$  & 0 \\
      6 & 2    & $[3_2,2]$ , $[3]$  &   0 \\
      6 & 2    & $[3_2]$      , $[3,2]$ & 0 \\
      6 & 2    & $[4,2_3]$ , $[2]$  & 0 \\
      6 & 2    & $[4,2_2]$ , $[2_2]$  &   0 \\
      6 & 2    & $[4]$      , $[2_4]$  &  $t+t^2$ \\
      6 & 3   & $[4]$ , $[2_3]$ , $[2]$ & 0 \\
      6 & 3   & $[4]$ , $[2_2]$ , $[2_2]$ & 0 \\
    \end{tabular}
\end{minipage}

\vspace{2mm}

For $d=7$ there is no complete classification (known by the authors).
The next table shows some examples (the first column provides either a possible
equation, or the reference where the corresponding curve has been constructed).

\vspace{2mm}

    \begin{tabular}{c|l|l|l|c}
           reference & $d$ & $\nu$   & type of cusp       & $-N(t)$ \\ \hline
   $x^6z+y^7$           & 7&    1        & $[6]$           &   0         \\
   $x^5z^2+y^7$         & 7&     2   & $[5,2_2]$ , $[2_3]$  &    0        \\
   $x^4z^3+ y^7$  &7&    2    & $[4,3]$ , $[3_2]$  & 0  \\
     \cite{Fenske} & 7&    2       & $[5]$ , $[2_5]$  & $2t+2t^2+2t^3$   \\
    \cite{Fenske}  & 7& 2         & $[4,2_3]$ , $[3_2]$  &  $t+2t^2+t^3$  \\
   \cite{Fenske} & 7&  2        & $[4,2_2]$ , $[3_2,2]$ & $2t+2t^2+2t^3$ \\
   \cite{Fenske}   & 7&2      &  $[4]$ , $[3_3]$    & $2t+2t^2+2t^3$   \\
   \cite{FZ2}  & 7&3     & $[4,2_2]$ , $[3_2]$ , $[2]$  & $3t+4t^2+3t^3$   \\
  \cite{FZ1} &7& 3 & $[5]$ , $[2_4]$ , $[2]$ & $3t+ 3t^3  $ \\
 \cite{FZ1} &7& 3 & $[5]$ , $[2_3]$ , $[2_2]$ & $3t+ 3t^3  $ \\
  \end{tabular}

\vspace{2mm}

\noindent
Here the first example is of Abhyankar-Moh-Suzuki type, while the next two
of Lin-Zaidenberg type,  cf. with the next sections.

\subsection{Example ($d$ large).}\label{dbig} Examples
with arbitrarily large $d$ and with non-vanishing $R(t)$ (but still satisfying
the conjecture) exist. E.g., if $C$ has  two cusps of types
$[d-2], [2_{d-2}]$ (see e.g. \cite{Fenske}), and $d$ is even, then
$$-N(t)=\sum_{k=1}^{\frac{d}{2}-2}\, k(t^k+t^{d-3-k}).$$
It is instructive (and sometimes rather mysterious) to verify  our conjecture
for the other  families listed e.g. in Fenske's article \cite{Fenske}.

\subsection{Remark.} 
Theorem 1 has the following immediate consequence:

\subsection{Corollary.}\label{cor} {\em Let $f=f_d+f_{d+1}+\cdots\ :
(\bc^3,0)\to (\bc,0)$  be a hypersurface superisolated singularity
with $\bar{\kappa}(\bp^2\setminus \{f_d=0\})<2$. Then the Seiberg-Witten
invariant  conjecture (cf. \cite{[51]}, see \ref{2.2} here) is true for
$(X,0)=(\{f=0\},0)$.}

\vspace{2mm}

\noindent It is a big question for the authors: what is the (conceptual)
connection between the log Kodaira dimension (of what?) with the general
conjecture  \ref{2.2} ?

\section{$\nu=1$ revisited. Comparison with other criterions.}

\subsection{}\label{nuegy}
Assume now that $\nu=1$, and write $Sing(C)=\{p\}$.
Recall that the characteristic polynomial $\Delta$ of $(C,p)\subset (\bp^2,p)$
is a complete (embedded) topological invariant of this germ, similarly as
the semigroup $\Gamma_{(C,p)}\subset \bn$, see \ref{a1}. In the next discussion
we will replace $\Delta $ by $\Gamma_{(C,p)}$.

Recall, that by \cite{CDG}, $\Delta(t)=(1-t)\cdot L(t)$, where $L(t)
=\sum_{k\in \Gamma_{(C,p)}} t^k$ is the
Poincar\'e series of $\Gamma_{(C,p)}$.
Since $P(t) $ was defined by the identity $\Delta(t)=1-P(t)(1-t)$, one gets
$L(t)+P(t)=1/(1-t)=\sum_{k\geq 0}t^k$. In particular,
$P(t)=\sum_{k \in \bn\setminus \Gamma_{(C,p)}} t^k$.

We let the reader to verify that Proposition \ref{abcde} implies the
following facts (cf. also with \ref{nuone}).

\subsubsection{Fact I}\label{fact1} {\em If $\nu=1$ then
$${\bf sw}(M)-\frac{K^2+s}{8}=\sum_{k\not\in\Gamma{(C,p)}}
\lceil k/d\rceil.$$
In particular, $p_g={\bf sw}(M)-(K^2+s)/8$ if and only if $D'(1)=0$,
or equivalently, if:}
$$\sum_{k\not\in \Gamma_{(C,0)}}\lceil k/d\rceil =d(d-1)(d-2)/6.$$

\subsubsection{Fact II}\label{fact2} {\em
$$Q(t)=\frac{(\sum_{k \not \in \Gamma_{(C,p)}} t^k )-\delta}{t-1}.$$
Then Conjecture \ref{2.4} (in the light of \ref{nuone})
predicts  that the coefficient $c_l$ of the monomial $t^{d(d-3-l)}$ of $Q(t)$
is exactly $(l+1)(l+2)/2$ (for any $0\leq l\leq d-3$). }

\subsubsection{Fact III}\label{fact3} {\em
$$D(t)=\sum_{k \in \Gamma_{(C,p)}} t^{\lceil k/d\rceil}-
\frac{1-t^d}{(1-t)^2}.$$
In particular, the second reformulation of the Conjecture \ref{2.4}
(and \ref{nuone}) is the identity
\begin{equation*}
\sum_{k \in \Gamma_{(C,p)}} t^{\lceil k/d\rceil}
=1+2t+\cdots+ (d-1)t^{d-2}+d(t^{d-1}+t^d+t^{d+1}+\cdots).\tag{$CP$}
\end{equation*}
}

\noindent Property ($CP$) (similarly as the property
of Fact II) connects the  local invariant $\Gamma_{(C,p)}$
with the degree $d$ of $C$.
It  predicts a very precise  distribution law for the elements of
the semigroup $\Gamma_{(C,p)}$ with respect to the intervals
$I_l:=(\,(l-1)d,ld\,]$ ($l\in \bn$). It says that for any $l\geq 0$ one has
\begin{equation*}\#\Gamma_{(C,p)}\cap I_l=\min\{l+1,d\}.\tag{$CP_l$}
\end{equation*}
By the symmetry of the semigroup (namely, $k\in \Gamma_{(C,p)}$
if and only if $2\delta-1-k\not\in \Gamma_{(C,p)}$), one has that
$(CP_l)$ is true if and only if $(CP_{d-2-l})$ is true.
In fact, $(CP_l)$ is automatically true for $l=0$ and any $l\geq d-2$.
But for $1\leq l
\leq  d-3$ it combines a lot of restrictions.

\subsection{A property of semigroups of plane curves}
\label{secgenus}

Now we will provide an alternative and very elementary proof of
Theorem~\ref{nuone}. It is based on a rather surprising and
general property about the semigroups of a (not necessarily
rational) plane curve at any collection of points with the sole
assumption that the curve is locally irreducible  at them.

Let $C$ be a (not necessarily rational) irreducible plane curve of degree $d$.
 Let $p_1,...,p_r$
be a set of points of $C$ such that the germ $(C,p_i)$ is irreducible for
any $i$.   Denote by $\Gamma_i$ the semigroup of $(C,p_i)$.

\begin{proposition}
\label{withgenus}
Suppose $l<d$. Let $(n_1,..., n_r)$ be
an r-uple of positive numbers such that $n_1+...+n_r\geq ld$. Then
\[\Sigma_{i=1}^r \#(\Gamma_i\cap [0,n_i])\geq (l+1)(l+2)/2.\]
\end{proposition}
\begin{proof}
Since the ideas of the proof appear when $r=1$, and this is the case what
we will use, for simplicity of notation we give
the proof only for that case.
Consider $n_1\geq ld$. Observe that $(l+1)(l+2)/2$ is the dimension of the
space of homogeneous polynomials $P$ in
three variables of degree $l$.
Let $\gamma(t)$
be a local parametrisation of $C$ at $p_1$. The composition $P(\gamma(t))$
can be written as
$P(\gamma(t))=\sum_{h=1}^\infty L_h(P) t^h$.
Therefore the conditions, imposed to $P$, that the
local intersection multiplicity $i_{p_1}(C,\{P\!=\!0\})$ is $>n_1$,
is given by the equations $L_h(P)=0$ for $0\leq h\leq n_1$.
On the other hand, as the semigroup $\Gamma_1$  is the collection of
intersection multiplicities of $C$ with other curves at $p_1$,
the number of
independent conditions is at most $\#(\Gamma_1\cap [0,n_1])$,
i.e. $L_h(P)=0$ for $h\leq n_1$ and $h\in \Gamma_1$
implies the vanishing $L_h(P)$ for all $h\leq n_1$.
Therefore, if
$\#(\Gamma_1\cap [0,n_1])<(l+1)(l+2)/2$, then there exists a
non-zero polynomial $P$ of degree $l$ with  $i_{p_1}(C,\{P\!=\!0\})>n_1$.
But, by Bezout's
theorem any such polynomial must have a component in common with $C$, but,
 as $C$ is irreducible, this is impossible.
\end{proof}

\begin{proof}[Proof of Theorem~\ref{nuone}]
Our conjecture~\ref{2.4} in the unicuspidal case implies that the polynomial
\[D(t)=\Sigma_{l=0}^{d-1} (\#\Gamma_1\cap ((l-1)d,ld]-(l+1))t^l\]
has the negative distribution property. We have proved in the previous proposition
that $D(t)$ has the positive distribution property. Hence
conjecture implies that $D(t)$ is identically $0$.
\end{proof}

\subsection{Comparing with consequences of log Bogomolov-Miyaoka-Yau type
inequalities.}\label{ap}  We are rather surprised
that the restrictions $(CP_l)$, imposed already by the very first intervals
$I_l$,
are closely related  with famous properties (conjectures) of rational
plane curves. In order to exemplify this, we need some more notations.

Similarly as above,
let $[m_1,m_2,\ldots]$ be the multiplicity sequence of $(C,p)$, the unique
singular point of $C$.
Following Nagata \cite{Na}, we define $t$ to be the maximal positive integer
such that $m_1\geq m_2+\cdots +m_t$.

We rewrite the  elements of $\Gamma_{(C,p)}$ as $\{0=\gamma_0<\gamma_1<
\gamma_2<\cdots\}$.
Recall that the minimal set of generators of
$\Gamma_{(C,p)}$ verifies $\bar{\beta}_0<\cdots < \bar{\beta}_g$.
In terms of the decorations of the Eisenbud-Neumann splice diagram  $\bar{\beta}_0=
p_1\cdots p_g=m_1$, $\bar{\beta}_k=a_kp_{k+1}\cdots p_g$ for $1\leq k<g$,
and $\bar{\beta}_g=a_g$, see \ref{a1}.
In fact, since $a_2>a_1p_1p_2$ and $p_1\geq 2$,
one has $\bar{\beta}_2>2\bar{\beta}_1$.
In particular, the first (at least) three
elements of $\Gamma_{(C,p)}$ depend only on
$\bar{\beta}_0$ and $\bar{\beta}_1$ (hence, all their properties can be
verified  essentially at the level of germs with one Puiseux pair).
Using these facts, one can verify that
\begin{equation*}
\left \{\begin{array}{l}
\gamma_1=m_1;  \ \ \gamma_2=m_1+m_2; \ \ \mbox{and}\\
m_1+2m_{t+1}\leq \gamma_3\leq 3m_1. \end{array}\right.\tag{5}\end{equation*}

\subsubsection{}\label{2.5}{\em Example. $(CP_1)$.} By the above notations,
 $(CP_1)$ says that $\gamma_2\leq d<\gamma_3$.

The first part
$\gamma_2\leq d$ can be verified easily. Indeed, if $L_p$ denotes the
tangent cone (line) of $C$ at $p$, then $i_p(C,L_p)\in \Gamma_{(C,p)}$ and
it is strict larger than $m_1$, hence
$\gamma_2\leq i_p(C,L_p)\leq d$ (the second inequality by B\'ezout's theorem).

For the second inequality, notice that if $C$ satisfies the Nagata-Noether
 inequality $d<m_1+2m_{t+1}$,
then by (5) it satisfies $d<\gamma_3$ too. Recall that any curve $C$
which can be transformed by a Cremona transformation into a line
satisfies the Nagata-Noether inequality \cite{Na,Ii,MS}. On the other hand,
it is conjectured that any cuspidal rational curve
can be transformed into a line  by a Cremona transformation (see \cite{MS}).

Notice also that via (5), $(CP_1)$ (or $d<\gamma_3$)
implies that $d<3m_1$, an inequality valid
for any cuspidal rational curve $C$, which was proved in \cite{MS} by
the log Miyaoka inequality \cite{My,KNS}. In fact, $d<\gamma_3$ also implies that
$i(C,L_p)=\gamma_2=m_1+m_2$.

\subsubsection{}\label{Yoscrit} Finally, let us reconsider again the inequality
$m_1+m_2=\gamma_2\leq d$. Here, by a result of Yoshihara \cite{Y1}
(whose proof is based on the Abhyankar-Moh theorem),
one has equality $m_1+m_2=d$ if and only if $C$ is an
Abhyankar-Moh-Suzuki curve (cf. \ref{3.1}).
The Abhyankar-Moh theorem \cite{AM} says that such a curve
can be transformed into a line, hence by the above discussion
it satisfies $(CP_1)$. In the next section we will show that, in fact,
it satisfies  all the restrictions $(CP_l)$, $l\geq 1$.

\subsubsection{}\label{Orevkovconj} The
interested reader is invited to  analyze some other
particular restrictions $(CP_l)$,  e.g. $(CP_l)$ for $l=1,2,3$, combined
together; cf. also with the last statement of \ref{4.3}(1).

\subsection{Remark.} Notice
also that the above distribution law $(CP)$ is very different from Arnold's
prediction of the distribution
  of {\em generic} sub-semigroups  of $\bn$ \cite{arnsem}.
This shows that the geometric realization of $(C,p)$ as the unique
singular point of a degree $d$ rational curve implies that the
semigroup of $(C,p)$ ``is far to be generic'' (a fact already
suggested by the local Abhyankar-Azevedo theorem as well).

\subsection{Comparison with Varchenko's criterion.}\label{varch} 
The above negative  distribution property has some  analogies with the
criterions provided by the semicontinuity of the spectrum
\cite{Var1,Var2}.  Namely, if $(C,p_i)_i$ are the local singular
points of the degree $d$ curve $C$, then the multisingularity
$\sum_i(C,p_i)$ is in the deformation of the universal plane germ
$(U,0):=(x^d+y^d,0)$. In particular, the collection of all
spectral numbers $Sp$ of the local plane curve singularities
$(C,p_i)_i$ satisfies the semicontinuity property compared with
the spectral numbers of $(U,0)$ for any interval
$(\alpha,\alpha+1)$.

In order to exemplify the similarities, let us assume
for simplicity that $\nu=1$.
 Since the spectral numbers of $(U,0)$ are of type
$l/d$, the semicontinuity property for intervals
$(-1+l/d,l/d)$ ( $l=2,3,\ldots , d-1$)  reads as follows:
\begin{equation*}
\#\{\alpha\in Sp\ :\ \alpha<l/d\}\leq (l-2)(l-1)/2.\tag{6}\end{equation*}
On the other hand,  the negative distribution property
of $D(t)$ is equivalent with
\begin{equation*}
\#\{k\in \Gamma_{(C,p)}\ :\ k\leq ld\}\leq
(l+1)(l+2)/2\tag{7}\end{equation*} for any $l=1,2,\ldots, d-3$.
Although the similarities are striking, it is not easy at all to
compare the two set of inequalities: the transition from
$Sp(C,p)$ to $\Gamma_{(C,p)}$ arithmetically is not very simple.
(If we add to the discussion the semicontinuity intervals
$(-1+l/d,l/d)$ with $l>d$, then the comparison is even more
difficult.) But even when $Sp(C,p)$ and  $\Gamma_{(C,p)}$ can be
easily compared (see e.g. below), still, the comparison of the
inequalities is not obvious. Let us exemplify this in the case
when $(C,p)$ has only one Puiseux pair, say $(a,b)$, $a<b$. In
this case the semicontinuity transforms  into
$$\#\{i\geq 0,\,j\geq 0\ :\ ia+jb<abl/d\}\leq (l-2)(l-1)/2 +[al/d]+[bl/d]+1.$$
(Recall that $(a-1)(b-1)=(d-1)(d-2)$.) On the other hand, our
criterion (7) is:
$$\#\{i\geq 0,\,j\geq 0\ :\ ia+jb\leq ld\}\leq (l+2)(l+1)/2.$$
In  the next examples show that the two restrictions are (at least
arithmetically) independent.

\subsection{Example.}\label{e2} {\em The semicontinuity of the spectrum
does not imply the negative distribution property.} Indeed,
  $(d,a,b)=(11,4,31)$ satisfies the semicontinuity but
 $N(t)$ has some positive coefficients
($\Gamma$ in $I_4$ has six elements). A similar example is
$(19,7,52)$; here  the first $l$ when $CP_l$ fails is $l=7$,
$I_l$ has 9 semigroup elements. (The reader is invited to verify
that these triples $(d,a,b)$ cannot be realized geometrically, a
fact compatible with our conjecture.)

In fact, one can even enter in the game another restriction, namely
the sharp  inequality of Orevkov \cite{Or}:
\begin{equation*}
d<\frac{3+\sqrt{5}}{2}(a+1)+\frac{1}{\sqrt{5}},\tag{8}\end{equation*}
and ask if the semicontinuity together with (8) would imply arithmetically
the negative
distribution property. The answer again is no, as (again)
the above two examples show.

Nevertheless, {\em for $d$ large,
$\nu=1$ and $(C,p)$ with only one Puiseux pair},
computer experiment shows that Orevkov inequality and the
semicontinuity of the spectrum imply the vanishing of $R(t)$.

\subsection{Example.}\label{537} Notice that for $\nu=1$,
\ref{nuone} replaces the negative distribution property with the
vanishing of $N(t)$. This is a  non-trivial additional
restriction. E.g.,   $(d,a,b)=(5,3,7)$ satisfies the
semicontinuity property of the spectrum and also  $N(t)=-t$. (Its
geometric realization is excluded by \ref{nuone}.)

\subsection{Example.}\label{e3} {\em The negative distribution property
does not imply the semicontinuity.}
Consider the case $(d,a,b)=(7,3,16)$.
Then it satisfies the negative distribution property (with $N(t)=-t^2$)
but it does not satisfy the semicontinuity of the spectrum.

Let us reconsider the triplet $(d,a,b)=(7,3,16)$ from the point of
view of \ref{nuone}: this
is a candidate for a unicuspidal curve but with $R(t)\not =0$.
According to \ref{nuone} a geometric realization cannot exist,
but the very existence of this triplet shows that the proof of \ref{nuone}
cannot be replaced by  a pure  arithmetical argument (having starting point
the non-positivity of the coefficients of $R$, and the fact that
$(C,p)$ has only one Puiseux pairs, and proving the vanishing of $R$).

\subsection{Remark.}\label{3.12} For $\nu=1$, we
do not know any example when $R$ is zero but the situation
geometrically is not realizable. (In particular, we don't know
any example when $R$ or $N$ are vanishing, but the semicontinuity
of the spectrum fails.)

{\em It would be a really remarkable characterization property, if
the conjecture \ref{fact2} (or, equivalently,
 \ref{fact3}), in fact,  would cover an
`if and only if' property.}

\section{The validity of $(CP)$ for Abhyankar-Moh-Suzuki curves.}

The \emph{AMS type} curves appeared naturally in the study of
rational plane curves $C$ meeting with a line $L$ at only one
point $\{p\}=C\cap L.$ By the Lin-Zaidenberg theorem \cite{ZL}
the curve $C$ has at most another cusp as singularities. $C$ is
called a {\em Abhyankar-Moh-Suzuki} curve if $C\setminus L$ is
smooth, respectively \emph{Lin-Zaidenberg} curve if $C\setminus L$
is singular.  K. Tono has given classifications, up to projective
equivalence, of {\em AMS} curves in \cite{To1} and of {\em LZ}
curves in \cite{To}. Here we will not use these classifications,
but we use the relation of these curves with automorphisms of the
affine plane $\bc^2=\bp^2\setminus L$ instead.

\subsection{Definition.}\label{3.1} An irreducible plane curve
$C$ is said to be of {\em Abhyankar-Moh-Suzuki type} ({\em AMS type} for
short) if there exists a line $L\subset \bp^2$ such that $C\setminus L
$ is isomorphic to $\bc$. In our situation this means that
 $\nu=1$ and $C\cap L=Sing(C)=\{p\}$.

Not any curve with $\nu=1$ is of {\em AMS} type, e.g.
the examples (c)-(f) in \ref{4.2} are not (cf. also with \ref{Yoscrit}).
The simplest {\em AMS} curve is $\{zx^{d-1}+y^d=0\}$.  In this case
 $\Gamma_{(C,p)}$ is generated by two elements, $d-1$ and $d$, and ($CP$)
can be easily verified. The goal of the present section is to prove
the general case.

\subsection{Theorem.}\label{3.2} {\em $(CP)$ is satisfied by
any AMS curve; in other words $R(t)\equiv 0$. }

\vspace{2mm}

Identify $\bc^2$ with $\bp^2\setminus L$. In the next discussion
we consider algebraic automorphisms   $\phi:\bc^2\to\bc ^2$
with components $(f,g)$. Recall, that by \cite{AM},
a curve $C$ is of {\em AMS} type if and only if  it is the compactification
in $\bp ^2$ of the zero locus of a component of a certain automorphism $\phi$.
The embedding of $\bc^2$ into $\bp^2$ allows us to view any automorphism
of $\bc^2$ as a birational transformation of $\bp^2$.
The point is that the combinatorics of the
minimal embedded resolution of $(C,p)$
is closely related to the combinatorics
of the minimal resolution of the indeterminacy of $\phi$ as a birational
transformation of $\bp^2$, and this last one can be described precisely.
In the sequel we make this statement more precise.
For details see ~\cite{Bo1,Bo2}.

Let $\pi':X'\to\bp ^2$ the minimal resolution of the indeterminacy of $\phi$.
Then  $(\pi')^*L$ is a divisor with normal crossings and its
dual graph $\calA$ (weighted by the corresponding self intersections)
has the following form:

\begin{center}
\includegraphics[width=12.5cm]{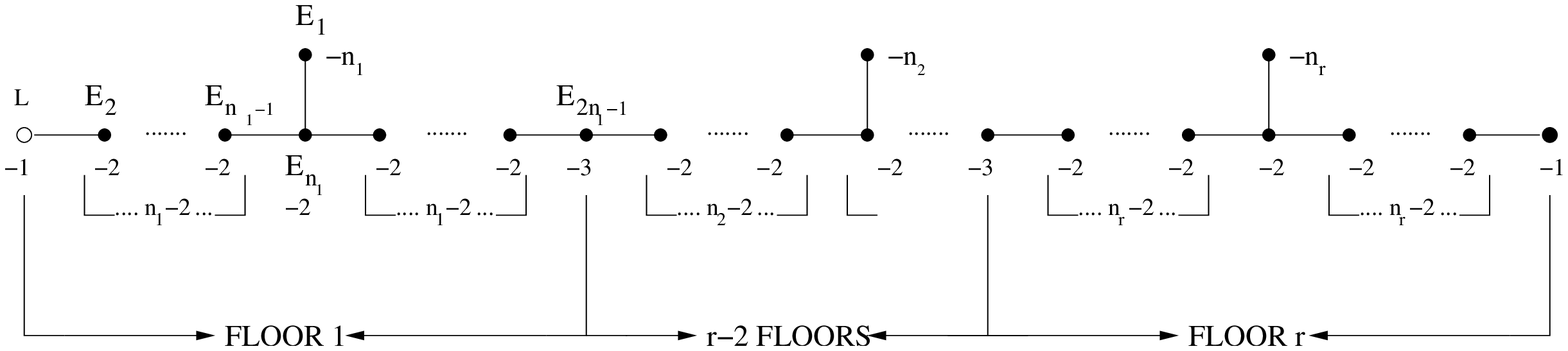}
\end{center}

\noindent
Here the vertices of $\calA$ are the ``black'' vertices.
The strict transform of $L$ is denoted by a ``white'' vertex with
decoration $L$. This graph can be obtained from elementary
pieces.

An {\em elementary graph of degree $n$} is by definition:

\begin{center}
\includegraphics[height=2.2cm]{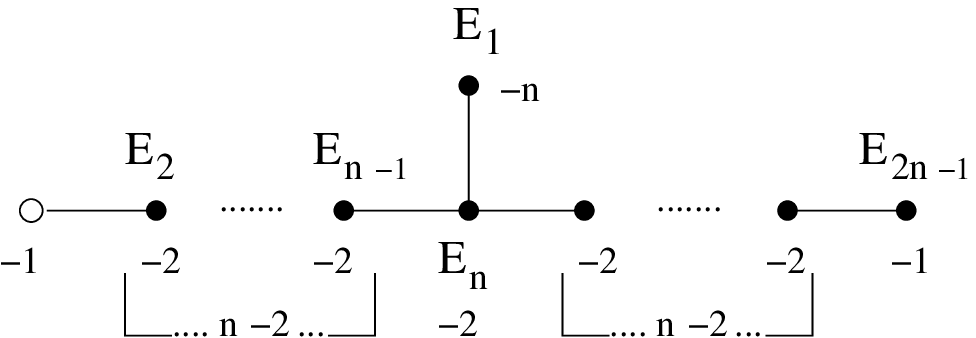}\label{figura}
\end{center}

\noindent
The graph $\calA$ is obtained by putting $r$ elementary graphs of degrees
$n_1,\ldots ,n_r$
one after the other, identifying the last vertex of each graph with
the first vertex of the next, and weighting the identified vertices with
$-3$. We say that the graph $\calA$ has $r$ floors. The morphism
$\pi':X'\to\bp^2$
is a composition of blowing ups, which are totally ordered by
appearance.  This gives a total order of the vertices of the
graph (which can be  recovered combinatorially by successively
contracting $-1$ vertices).

Now, we concentrate on the unique singular point $p$ of $C$.
The minimal embedded resolution  $\pi:X\to\bp^2$ of $(C,p)\subset (\bp^2,p)$
is again a composition of blowing ups, and  all these
blowing ups appear in the minimal resolution of the indeterminacy of $\phi$.
In fact, there are two possibilities:

\vspace{1mm}

\no
(a) \ $\pi$ is the composition
of all the blowing ups of $\pi'$ except the last $n_r-1$ of them;\\
(b) \ $\pi$ is the composition of the blowing ups of the first $r-1$
floors of $\pi'$ except the last  $n_{r-1}-1$ of \hspace*{5mm} them.

\vspace{1mm}

\noindent
If the second possibility holds then one can find a different automorphism
$\psi=(f,g)$ such that $C$ is the closure of $\{f=0\}$ and the graph
associated with $\psi$ has $r-1$ floors satisfying possibility (a).
Therefore, one can always assume the validity of (a).
This means that the embedded resolution graph of $(C,p)$ can be obtained from
the above graph $\calA$ by deleting the last $n_r-1$ (black) vertices, and
changing the decoration of $E_{n_r}$ into $-1$. The strict transform
$\bar{C}$ intersects $E_{n_r}$.

In particular, the singularity $p$ of the curve $C$ has $r$ Newton pairs,
and the decorations
$\{(p_k,a_k)\}_{k=1}^r$ of its
Eisenbud-Neumann splice diagram (cf. \ref{a1}) can be computed  by standard graph-determinant
computations (cf. \cite{EN}, section 21). One gets the following:
\begin{equation*}\left\{ \begin{array}{l}
d=i_p(z,\bar{C})=n_1\cdots n_r;\\
p_1=n_1-1;\ \mbox{and} \ p_k=n_k \ \mbox{for $k=2,\ldots,r$};\\
a_1=n_1, \ \mbox{and} \ a_k=n_1^2\cdots n_{k-1}^2n_k-1 \ \mbox{for $k
=2,\ldots , r$}.\end{array}\right.\tag{9} \end{equation*}
In particular, the generators $\bar{\beta}_i$ of $\Gamma_{(C,p)}$ are
$\bar{\beta}_0=(n_1-1)n_2\cdots n_r$, $\bar{\beta}_1=n_1\cdots n_r$, and
for $1<k\leq r$  one has $\bar{\beta}_k=(n_1^2\cdots n_{k-1}^2n_k-1)n_{k+1}
\cdots n_r$.

We will verify $(CP)$ by induction over $r$.
Assume that $\calA=\calA_r$ has $r$ floors
of degrees $n_1,\ldots,n_r$. Let $\Gamma_r$ be the corresponding
semigroup. If $r=1$ then $\Gamma_{1}$ is generated by $n_1-1$
and $n_1$, $d_1=n_1$ and $(CP)$ can be easily verified. If $r>1$,
by (9) and \ref{a1}, the degree $d_r$ of
of $C$ is $d_r=n_1\cdots n_r=d_{r-1}n_r$, and
$\Gamma_{r}$ is generated by
$n_r\Gamma_{r-1}$ and $\bar{\beta}_r:=n_1^2\cdots n_{r-1}^2n_r-1$.
In fact, any element $x\in \Gamma_r$ can be written in a unique way
as $x=n_r y+b\bar{\beta}_r$ with $y\in \Gamma_{r-1} $ and $0\leq b< n_r$
(see e.g. \cite{TZ}).

Now, the inductive step runs as follows:
$$\sum_{x\in\Gamma_r}t^{\lceil x/d_r\rceil}=
\sum _{b=0}^{n_r-1}\sum_{y\in \Gamma_{r-1}}
t^{\lceil (yn_r+b\bar{\beta}_r)/d_r\rceil}=
\sum _{b=0}^{n_r-1}\sum_{y\in \Gamma_{r-1} }
t^{bd_{r-1}+\lceil y/d_{r-1}\rceil}=
\frac{1-t^{d_r}}{1-t^{d_{r-1}}}\cdot \frac{1-t^{d_{r-1}}}{(1-t)^2}.$$

\section{The case $\nu=1$ and  $(C,p)$ with one Puiseux  pair}

\subsection{}\label{4.1} Assume that the unique cusp of
$C$ has exactly one Puiseux  pair $(a,b)$ (where $1<a<b$).
Then clearly $(a-1)(b-1)=(d-1)(d-2)$, where $d=\deg(C)$ as above.

In the sequel we denote by $\{\varphi_j\}_{j\geq 0}$ the Fibonacci numbers
$\vp_0=0$, $\vp_1=1$, $\vp_{j+2}=\vp_{j+1}+\vp_j$.
The Fibonacci numbers have a remarkable amount of interesting properties,
see e.g. Vajda's book \cite{Vaj}.
We also will use some of them which will be crucial in the next arguments
and also in the proof of the conjecture for Kashiwara's families
(cf. section 6). E.g.:

\vspace{1mm}

(i) $\vp_j^2-\vp_{j-1} \vp_{j+1}=(-1)^{j+1}$;

(ii) $\vp_j^2-\vp_{j-2} \vp_{j+2}=(-1)^{j}$.

(iii) $\gcd(\vp_j,\vp_i)=\vp_{\gcd(j,i)}$.

\subsection{Examples.}\label{4.2} We will consider the following pairs
$(a,b)$:

\vspace{1mm}

(a) $(a,b)=(d-1,d)$;

(b) $(a,b)=(d/2, 2d-1)$, where $d=\deg(C)$ is even;

(c) $(a,b)=(\vp_{j-2}^2,\vp_j^2)$ and $d=\vp_{j-1}^2+1=\vp_{j-2} \vp_{j}$, where $j$ is odd
and $\geq 5$;

(d) $(a,b)=(\vp_{j-2},\vp_{j+2})$ and $d=\vp_{j}$, where $j$ is odd and
$\geq 5$;

(e) $(a,b)=(\vp_4,\vp_8+1)=(3,22)$ and $d=\vp_6=8$;

(f) $(a,b)=(2\vp_4,2\vp_8+1)=(6,43)$ and $d=2\vp_6=16$.

\vspace{1mm}

All these cases are realizable: (a) e.g. by $\{zy^{d-1}=x^d\}$,
(b) by $\{(zy-x^2)^{d/2}=xy^{d-1}\}$, or by the parametrization
$[z(t):x(t):y(t)]=[1+t^{d-1}:t^{d/2}:t^d]$. Both cases (a) and (b) satisfy
Yoshihara's criterion $m_1+m_2=d$ \cite{Y1} (cf. here with \ref{Yoscrit}),
hence any curve with these data is {\em AMS} curve.

The existence of (c) and (d) is guaranteed by Kashiwara classification
\cite{Kash}, Corollary 11.4. These two cases can be realized by a rational
pencil of type $(0,1)$: the generic member of the pencil is (c), while the
special member of the pencil is of type (d) (cf. also with the next section).
Orevkov in \cite{Or}
denoted the curves (d) by $C_j$, where a different construction is also given
for them.

The cases (e) and (f) correspond to  the sporadic cases $C_4$ and $C_4^*$ of
Orevkov \cite{Or}.

\subsection{Remarks.}\label{4.3} (1) The above list is not accidental.
In \cite{class} we prove that if $C$ is a unicuspidal rational plane curve
of degree $d$ and if the singular point $p$ of $C$ has only one
characteristic pair $(a,b)$,
then the triple $(d,a,b)$ is one of the above cases.
This classification is coordinated by the following integer.
Let $\pi:X\to \bp^2$ be the minimal embedded resolution of $C\subset
\bp^2$, and let $\bar{C}$ be the strict transform of $C$.
Then $\bar{C}^2=d^2-ab=3d-a-b-1$, and in
the above cases is as follows: it is positive for (a) and (b),
it is zero for (c), equals $-1$ for (d), and $=-2$ for (e) and (f).
[Notice that $\bar{C}^2<-1$ if and only if $a+b>3d$, i.e. if and only if
the semigroup element $a+b$ is not sitting in the first three intervals
$I_l$. Notice also that $\bar{C}^2<-1$ happens exactly for the cases (a)-(d),
i.e. when
$\bar{\kappa}(\bp^2\setminus C)=-\infty$; cf. (3) below.]

(2) Since $\Gamma_{(C,p)}$ is the semigroup generated by $a$ and $b$,
the verification of $(CP)$ for the above triples $(d,a,b)$
is purely combinatorial depending on these integers.

On the other hand, not any triple $(d,a,b)$ (with $(a-1)(b-1)=(d-1)(d-2)$)
satisfies  $(CP)$. E.g., $(5,3,7)$ or  $(17,6,49)$ do not. (But
curves with these data do not exist, cf. \cite{class}.)

(3) The
log Kodaira dimensions  $\bar{\kappa}(\bp^2\setminus C)$ are the following
(cf. \cite{Or}): $-\infty$ for (a)-(d), and 2 for the last two
sporadic cases.

(4) Let $\alpha=(3+\sqrt{5})/2$
be the root of $\alpha+\frac{1}{\alpha}=3$. Notice that in family (d)
$d/a$ and $b/d$ asymptotically equals $\alpha$.
In fact, for $j$ odd, $\{\vp_j/\vp_{j-2}\}_j $ are the increasing convergents
of the continued fraction of $\alpha$. Using this,
another remarkable property of the family (d) can be described as follows
(cf. \cite{Or}, page 658).
The convex hull of all the pairs $(m,d)\in \bz^2$
satisfying $m+1\leq d<\alpha m$ (cf. with the sharp Orevkov inequality
\cite{Or}, or \ref{e2})  coincides  with the convex hull of all pairs $(m,d)$
realizable by rational unicuspidal curves $C$ (where $d=\deg(C)$ and
$m=mult(C,p)$) with
$\bar{\kappa}(\bp^2\setminus C)=-\infty$; moreover, this convex hull is
generated by curves with numerical data (a) and (d). (For curves with
$d>\alpha m$, see \ref{5.1}.)

\subsection{Theorem.}\label{4.4} {\em The identity $(CP)$ (i.e. $R(t)\equiv 0$)
is satisfied in the above  cases (a)-(f).}

\vspace{2mm}

The cases (a) and (b) are covered by \ref{3.2}. The cases (e) and (f)
can be verified by hand: doing this the reader definitely will feel
the mystery of this distribution pattern. In the sequel we will verify
(c) and (d).

We start with (d). We fix $j$. The point is that $a+b=3d$,
hence a nice induction can be considered if we group the intervals $I_l$
in blocks of three. Let us analyze the first block.
$I_1$ contains $a$ and $2a$; $I_2$ contains $3a$, $4a$, and $5a$; finally
$I_3$ has $6a$, $7a$, $b$ and $a+b$. All this can be verified by the
definition of the Fibonacci numbers.
E.g. $3a>d$  iff $3\varphi_{j-2}>\varphi_{j-1}+\varphi_{j-2}$
iff $2\varphi_{j-2}>\varphi_{j-1}$ iff $2\varphi_{j-2}>\varphi_{j-2}+
\vp_{j-3}$ which is true.

It is clear that it is enough to analyze the intervals $I_l$ for $l<d$.
So fix such an $l$.
The point is that an inequality of type $kb>ld$ is true if and only if
$k/l>1/\alpha$. One direction is easy, since $d/b>1/\alpha$. Assume, that
$k/l>1/\alpha$, then $d/b\geq k/l>1/\alpha$ is not possible since
$d/b$ -- being a convergent of $\alpha$ -- is the best approximation of
$\alpha$ among fractions with denominator $\leq b$
(see e.g. \cite{Niven} (7.13)). Hence $k/l>d/b$.
In particular, $kb\in I_l$ if and only if $\lceil k\alpha\rceil =l$.

There is a similar statement for $ia$, but its proof is slightly different.
If $ia\leq ld$ then using $d/a<\alpha$ one gets $i/l<\alpha$.
Assume now that $i/l<\alpha$ then we wish to prove that $d/a<i/l<\alpha$
is not possible. For this notice that $d/a<b/d<\alpha$ (in fact $d/a$ and
$b/d$ are two consecutive convergents of $\alpha$). Then $i/l$ cannot be
between $b/d$ and $\alpha$ since $b/d$ is one of the convergents and $l<d$;
and also $i/l$ cannot be between $d/a$ and $b/d$ since $ab-d^2=1$ and
$l<d$ (cf. \cite{Niven}, p. 165, 5(a)). Hence $ia\leq ld$ if and only if
$i/l<\alpha$. This implies that $ia\in I_l$ if and only if
$\lceil i/\alpha\rceil =l$.

Now, using the fact that $a+b=3d$, we can move the first block to the
second block (of three intervals ) by adding $a+b=3d$.
In the relevant intervals ($l<d$)
the only terms not in the image of this translation have form
$ia$ or $kb$, and they are all distinct.
An easy counting shows that the only fact we have to show is
that in any interval there are exactly 3 terms of these types.
Namely, we have to verify the following:
Consider the numbers  $S_1$ of the form $\lceil i/\alpha\rceil$
and the numbers $S_2$ of the form  $\lceil k \alpha\rceil$.
Then the claim is that in $S_1\cup S_2$
each positive integer  appears exactly three times.

This can be proved as follows. Consider in the $(x,y)$-plane the semi-line
$\ell_1:$
$\{y=\alpha x\}$ (in the first quadrant) and the semi-line $\ell_2$:
$\{y=-x/\alpha\}$ (in the forth quadrant).
For any positive integer $l$ consider the vertical segment
connecting the two intersection points of the line $\{x=l\}$ with
$\ell_1 $ and $\ell_2$. Notice that the length of this segment is
exactly $3l$. Any horizontal line $\{y=i\}$ (resp. $\{y=-k\}$)
 intersects the segment iff $i/\alpha\leq l$ (resp.
$k\alpha \leq l$). Since the segment can be intersected exactly by
$3l$ horizontal lines of type $y=$integer (and all the numbers of type
$k\alpha$ and $i/\alpha$ are distinct) the claim follows.

\vspace{2mm}

Next, we prove (c). We start as in (d).
Fix $j$ and notice that $d=\vp_{j-2}\vp_j$.

Step 1. Assume that $l<\vp_{j-2}+\vp_j$. If $kb>ld$ then $k\vp_j>l\vp_{j-2}$,
hence $k/l>\vp_{j-2}/\vp_j>1/\alpha$. Conversely, assume that
$k/l>1/\alpha$. Consider the convergents
$\vp_{j-2}/\vp_j>\vp_j/\vp_{j+2}>1/\alpha$. Since $\vp_{j+2}>l$
we conclude that $k/l$ cannot be neither
between $\vp_j/\vp_{j+2}$ and $1/\alpha$
nor between $\vp_{j-2}/\vp_j$ and $\vp_j/\vp_{j+2}$
by similar argument as in (d). Hence either $k/l>\vp_{j-2}/\vp_j$ or
$l=\vp_j$. This shows that for $l\not\in \{ \vp_j,\vp_j+1\}$,
$kb\in I_l$ if and only if $\lceil k\alpha\rceil =l$.

Step 2. Assume again that $l<\vp_{j-2}+\vp_j$. If $ia\leq ld$
then $i/l\leq \vp_j/\vp_{j-2}<\alpha$. Conversely, assume that
$i/l<\alpha$ and consider the three intervals
$\vp_j/\vp_{j-2}<\vp_{j+2}/\vp_j<\vp_{j+4}/\vp_{j+2}<\alpha$.
Then $i/l$ cannot be in the second and third interval by similar arguments as
above. Moreover, it cannot be in the first one either,
since the two end-points
are two elements of the $\vp_{j}$-Farey sequence, and the denominator
of any rational number between them
is at least $\vp_{j-2}+\vp_j$ (cf. (6.4) \cite{Niven}).
Therefore, either $l=\vp_j$ or $i/l\leq \vp_j/\vp_{j-2}$.
Hence, if $l\not\in \{\vp_j,\vp_j+1\}$ then
 $ia\in I_l$ if and only if $\lceil i/\alpha \rceil=l$.

Step 3. The intervals $I_l$ for $l=\vp_j$ and $\vp_j+1$ can be analyzed
independently: $kb\in I_{\vp_j}$ if and only if $k=\vp_{j-2}$
(and then $kb=\vp_j d$), and $ia\in I_{\vp_j}$ if and only if
$i=\vp_{j+2}-1$ or $i=\vp_{j+2}-2$.

Similarly, $I_{\vp_j+1}$ contains no number of type $kb$, but contains
three numbers of type $ia$ for $i=\vp_{j+2},\vp_{j+2}+1,\vp_{j+2}+2$.

Combined the argument of (d) applied for Step 1 and Step 2
and the above two facts, we  conclude that the distribution
patern is true for any $l<\vp_j+\vp_{j-2}$.

Step 4. Notice that in the relevant intervals (i.e. for $l<d$),
the semigroup element $ia+kb$ determines uniquely the integers $i$ and $k$.

Step 5. Using the identity $a\vp_j=\vp_{j-2}d$, we construct a well-defined
injective map $s_l:I_{l-\vp_{j-2}}\cap \Gamma \to I_l\cap \Gamma$
by $x\mapsto x+a\vp_j$. Denote by $P_l$ the subset of semigroup elements
of $I_l\cap \Gamma$  which have the form $ia+kb$ with $k\geq 0$ and $
a>i\geq 0$. Then (via step 4)
$I_l\cap \Gamma$ is the disjoint union of the image of $s_l$ and $P_l$.
Therefore, it is enough to show that
\begin{equation*}
\#P_l=\vp_{j-2}\tag{10}
\end{equation*}
for any $\vp_{j-2}\leq l<d$.

Step 6. Since the distribution property is true for any
$l<\vp_{j-2}+\vp_j$ (step 3), we conclude (by the arguments of step 5)
that (10) is true for any $\vp_{j-2}-1\leq l\leq \vp_{j-2}+\vp_j-1$.

Step 7. We verify that there is a bijection $P_l\to P_{l+\vp_j}$
given by $(i,k)\mapsto (i,k+\vp_{j-2})$ for any $l\geq \vp_{j-2}$.
The facts that the map is well-defined and injective are clear, for the
surjectivity one has to verify an (easy)
inequality satisfied by the Fibonacci  numbers.
Hence (10) follows by step 6 and induction.

\section{$R(t)\equiv 0$ for
$\bar{\kappa}(\bp^2\setminus C)=-\infty$ (i.e. for Kashiwara's curves)}

We start with the following lemma which helps  us to
compare distinct semigroups.

\subsection{Lemma.}\label{compare} {\em Fix integers $k\geq 0$ and $m,d\geq 1.$ Then

(a) in the set
$B:=\{\lceil m k/d\rceil+j\}_{j=0}^{m-1}$ the unique multiple of $m$ is
${\lceil k/d\rceil} m.$

(b) Fix a subset  $\Gamma\subset \bn$ and define the series
\begin{equation*}
\eta(t):=\sum _{k\in\Gamma}\, t^{\lceil k/d\rceil} \, \, \mbox{and } \,\,
\chi(t):=\sum _{k\in\Gamma}\, t^{\lceil m k/d\rceil}.
\end{equation*}
Then  the series $\psi(t):=\chi(t)(1-t^m)/(1-t)$ satisfies the identity}
\begin{equation*}
\eta(t^m)=\frac{1}{m} \sum _{\xi^m=1} \psi(\xi t).
\end{equation*}
\begin{proof}
It is clear that there is only one multiple of $m$  in $B,$ say $a$. If $k=qd$ the proof is also clear. Otherwise, write
$k=qd+r$ with $r\in\{1,\ldots,d-1\}.$ Thus $mq/d=mq+rm/d.$ Since $0<r/d<1$ then $0<rm/d<m$ and
$0< \lceil m r/d\rceil \leq m.$  Finally $mq<a\leq mq+m+m-1<m(q+2),$ and $a=(q+1).$
\end{proof}

\subsection{The proof of the conjecture for Kashiwara's curves.}\label{6.2} \
As we have mentioned most of the previous curves appear
as irreducible components of fibres of rational functions on $\bp^2$
of type $(0,1)$, that is, rational functions all whose fibres (once the indeterminacy point has been
removed) are isomorphic to $\bc.$ Remark that any of these rational function $\phi$
has at most one indeterminacy point which will be the only possible singular point of any of its fibres.
The classification of such rational functions is given in H. Kashiwara's paper
\cite{Kash}. Moreover, a rational cuspidal curve $C$
verifies $\bar{\kappa}(\bp^2\setminus C)=-\infty$ if and only if $C$ is an irreducible component
of a fibre of a rational function of type $(0,1).$

The number of multiple fibres that a rational function $\phi$ on $\bp^2$
of type $(0,1)$ can have is at most two.
Kashiwara's classification gives three strata inside the set of rational functions of type $(0,1)$:
${\mathcal F}_0, {\mathcal F}_I, {\mathcal F}_{II},$ according to the number of multiple fibres will be $0, 1$ or $2$
respectively.

\vspace{1mm}

(1) The stratum ${\mathcal F}_0$ consists of all linear rational functions
on $\bp^2$, therefore there are no cuspidal rational curves as fibres.

(2) If $\phi\in {\mathcal F}_I$ then the multiple fibre of $\phi$ turns out to be a line $L$, and $\phi$ is a component of an automorphism of $\bc^2=\bp^2\setminus L$, see Corollaire 8.1 in \cite{Kash}. Moreover every fibre of $\phi$, but the multiple one, defines a rational
cuspidal curve of {\emph AMS-type}
for which we have already checked $(CP)$ in \ref{3.2}.

(3) A rational function $\phi\in {\mathcal F}_{II}$ has two atypical  values, say $0$ and $\infty$.
Its divisor will be denoted by $(\phi)=mS_0-nS_\infty$ (for some integers that we may assume $\gcd(m,n)=1$ taking $\phi$ primitive). From the topological point of view the rational function $\phi$
has only three different fibres $S_0$, $S_\infty$ and the \emph{generic fibre} $\phi_{ge}.$

The minimal resolution graph of the indeterminacy of a rational function $\phi$ of type $(0,1)$ almost coincides with the minimal resolution graph of the generic fibre of $\phi$. But this is not the case for special fibres.
In her Th\'eor\`eme 6.1 the minimal resolution graphs of the indeterminacy point
of the rational functions of type $(0,1)$  are given. The remaining case in which we are interested in is   $\phi\in {\mathcal F}_{II}$ where
five different graphs appear:

Case 1. $II(\ell)^*$ (with $\ell\geq 0),$

Case 2. $II^+(\ell,N;\lambda_1,\ldots,\lambda_N)^*$ (with $N$ even $\geq 2),$

Case 3. $II^+(\ell,N;\lambda_1,\ldots,\lambda_N)^*$ (with $N$ odd $\geq 1),$

Case 4. $II^-(\ell,N;\lambda_1,\ldots,\lambda_N)^*$ (with $N$ even $\geq 2),$

Case 5. $II^-(\ell,N;\lambda_1,\ldots,\lambda_N)^*$ (with $N$ odd $\geq 2).$

Here $\lambda_1,\ldots,\lambda_N$ are integers such that $\lambda_1,\ldots,\lambda_N\geq 0$
if $\ell\geq 1$ and $\lambda_1,\ldots,\lambda_N\geq 1$
if $\ell=0.$

In Case 1, for a rational function $\phi^\ell$ whose graph belongs to $ II(\ell)^*,$ there are three
distinct cuspidal rational curves appearing as fibres of $\phi^\ell$:
$S_0^\ell$, $S_\infty^\ell$ and the generic fibre $\phi_{ge}^\ell$.
The rational function $\phi^\ell$ is constructed by induction based on a rational function $\phi^{\ell-1}$ which implies that
$S_\infty^\ell$ is nothing but $S_0^{\ell-1}$, see e.g. Corollaire 11.4.
Thus for each graph in $ II(\ell)^*$ there are only two new unicuspidal rational plane curves:
the generic member of the pencil, we will denote it by $II(\ell)_{ge},$
and the new special member of the pencil, we will denote it by $II(\ell)_{sp}.$
In examples (c) and (d) in \ref{4.2} we have already studied all these curves:

(c) $II(\ell)_{ge}$ has degree $d=\vpt \vpf$ and only one characteristic pair $(a,b)=(\vpt^2,\vpf^2);$

(d) $II(\ell)_{sp}$ has degree $d=\vpt$ and only one characteristic pair $(a,b)=(\vpo,\vpf).$

Therefore in \ref{4.4} we have proved the identity $(CP)$ for $II(\ell)_{ge}$
and $II(\ell)_{sp},$ $\ell\geq 0$.

The same fact happens for all rational functions whose graphs belongs to any of
the remaining four cases, see e.g.
Corollaire 11.6. The irreducible component of the $\infty$ fibre of the
${\ell}$-rational function is nothing but
the $0$-fibre of the previous ${\ell-1}$-rational function.
Therefore in each group we have just two new unicuspidal rational curves:
the generic member of the pencil and one special member of the pencil.

The only difference between the generic
rational curve in  $II^+(\ell,N;\lambda_1,\ldots,\lambda_N)^*$
(with $N$ even $\geq 2),$ (Case 2), and the generic rational curve in
$II^+(\ell,N;\lambda_1,\ldots,\lambda_N)^*$ (with $N$ odd $\geq 1),$ (Case 3),
is the number of Puiseux pairs.
In fact, it is possible to codify the invariants of the generic  members of Cases 2 and 3 in one and the same sequence of Eisenbud-Neumann splice diagrams, which we denote by  $II^+(\ell,N;\lambda_1,\ldots,\lambda_N)_{ge}.$
The same can be done with the special members of the pencils, which we will denote by
$II^+(\ell,N;\lambda_1,\ldots,\lambda_N)_{sp}.$

Finally the same can be done also with Cases 4 and 5, that is for rational functions
with graphs in $II^-(\ell,N;\lambda_1,\ldots,\lambda_N)^*$ (with $N$ even $\geq 2)$  or $II^-(\ell,N;\lambda_1,\ldots,\lambda_N)^*$ (with $N$ odd $\geq 2).$
The generic (resp. special) curve will be denoted by
$II^-(\ell,N;\lambda_1,\ldots,\lambda_N)_{ge}$ (resp.   $II^-(\ell,N;\lambda_1,\ldots,\lambda_N)_{sp}$).

\subsection{Degrees of curves and generators of semigroups}
Fix $\ell\geq 0$ and $N\geq 1$ and  non-negative integers $\lambda_1,\ldots,\lambda_N$.

\subsubsection{$II^+(\ell,N;\lambda_1,\ldots,\lambda_N)_{ge}$}
For $1\leq i\leq N$ define

\begin{equation*}\left\{ \begin{array}{l}
n_i:=\lambda_i \vpt^2 +\vpt\vpm-1  \qquad \mbox{for $i$ odd};\\
n_i:=\lambda_i \vpt^2 +\vpt(\vpt-\vpm)-1  \ \mbox{for $i$ even}.
\end{array}\right.\tag{11} \end{equation*}

The degree of the curve $II^+(\ell,N;\lambda_1,\ldots,\lambda_N)_{ge}$ is
$d=\vpt\vpf n_1\cdots n_N,$ (see Proposition 7.2 in \cite{Kash}). We remark that the integer $m_\ell$ in \cite{Kash} is nothing but the Fibonacci number $\vpt.$
This curve has only one singularity  whose splice diagram
can be easily deduce from the resolution graph of the corresponding rational function.
This singularity has $N+1$ Newton pairs and the decorations $\{(p_k,a_k)\}_{k=1}^{N+1}$
of the corresponding Eisenbud-Neumann splice diagram can be done computing the graph determinants (as it is explained in \cite{EN}).
Thus

\begin{equation*}\left\{ \begin{array}{l}
p_k=n_k \ \mbox{for $1\leq k \leq N,$} \,\mbox{and} \ p_{N+1}=\vpt^2;\\
a_k=(\vpf^2 n_1^2\cdots n_{k-1}^2n_k-1)/{\vpt^2}, \ \mbox{for $1\leq k \leq N$}\ \mbox{and} \ a_{N+1}=\vpf^2 n_1^2\cdots n_{N}^2.
\end{array}\right.\tag{12} \end{equation*}

To check that the $a_k$'s are integers, let $b_k$ denote the numerator $\vpf^2 n_1^2\cdots n_{k-1}^2n_k-1$ of $a_k,$ then $b_k=b_{k-1} n_{k-1}n_k+n_{k-1}n_k-1.$
By induction, to check that $b_k= 0 \mod \vpt^2$
it is enough to prove that $n_{k-1}n_k-1= 0 \mod \vpt^2.$ From the definition of $n_k$, see $(11)$,
$n_{k-1}n_k-1=(\vpt\vpm-1)(-\vpt\vpm-1)-1=-(\vpt^2\vpm^2-1)-1=0\, \mod \vpt^2,$ after property $(ii)$ of
the Fibonacci numbers.

In particular, the generators $\{\bar{\beta}_k\}$ of its semigroup
$\Gamma_{II^+(\ell,N;\blambda)_{ge}}$ are
$\bar{\beta}_0=\vpt^2 n_1\cdots n_N$, $\bar{\beta}_k=(\vpf^2 n_1^2\cdots n_{k-1}^2n_k-1)n_{k+1}\cdots n_N$
for $1\leq k\leq N$,
and $\bar{\beta}_{N+1}=\vpf^2 n_1^2\cdots n_{N}^2.$

\subsubsection{$II^+(\ell,N;\lambda_1,\ldots,\lambda_N)_{sp}$}

The unicuspidal rational curve $II^+(\ell,N;\lambda_1,\ldots,\lambda_N)_{sp}$ has degree
$d=\vpf n_1\cdots n_N,$ (see Proposition 7.2 in \cite{Kash}).
This curve has also only one singularity which is a cusp with $N+1$ Newton pairs. To compute the decorations
$\{(p_k,a_k)\}_{k=1}^{N+1}$of the corresponding splice diagram
all the determinants, but the last two, are the same as in the previous case.
Thus the invariants in this case are:

\begin{equation*}\left\{ \begin{array}{l}
p_k=n_k \ \mbox{for $1\leq k\leq N,$} \ \mbox{and} \ p_{N+1}=\vpt;\\
a_k=(\vpf^2 n_1^2\cdots n_{k-1}^2n_k-1)/{\vpt^2},  \mbox{for $1\leq k\leq N$}\, \mbox{and} \ a_{N+1}=(\vpf^2 n_1^2\cdots n_{N}^2+1)/{\vpt}.
\end{array}\right.\tag{13} \end{equation*}

Moreover $a_{N+1}$ is also integer because
$n_i=-1 \mod \vpt$ and then $\vpf^2 n_1^2\cdots n_{N}^2+1=\vpf^2-1=(\vpt\vps-1)+1=0 \mod \vpt.$

The generators $\{\bar{\beta}_k\}$ of its semigroup $\Gamma_{II^+(\ell,N;\blambda)_{sp}}$ are
$\bar{\beta}_0=\vpt n_1\cdots n_N$, $\bar{\beta}_k=(\vpf^2 n_1^2\cdots n_{k-1}^2n_k-1)n_{k+1}\cdots n_N/\vpt$
for $1\leq k\leq N$,
and $\bar{\beta}_{N+1}=(\vpf^2 n_1^2\cdots n_{N}^2+1)/\vpt.$

\subsubsection{$II^-(\ell,N;\lambda_1,\ldots,\lambda_N)_{ge}$}
For $1\leq i\leq N,$  define

\begin{equation*}\left\{ \begin{array}{l}
{\tilde n}_i:=\lambda_i \vpt^2 +\vpt\vpm-1  \qquad \mbox{for $i$ even};\\
{\tilde n}_i:=\lambda_i \vpt^2 +\vpt(\vpt-\vpm)-1  \ \mbox{for $i$ odd}.
\end{array}\right.\tag{14} \end{equation*}

\no
Essentially the data of the curve $II^-(\ell,N;\lambda_1,\ldots,\lambda_N)_{ge}$ can be obtained from the data of
the curve $II^+(\ell,N;\lambda_1,\ldots,\lambda_N)_{ge}$ replacing $\vpf$ by $\vpo.$ The combinatorial reason for that
is the identity $\vpt^2=\vpo \vpf -1.$
The degree of the rational curve $II^-(\ell,N;\lambda_1,\ldots,\lambda_N)_{ge}$ is
$d=\vpt\vpo {\tilde n}_1\cdots  {\tilde n}_N,$ (see Proposition 7.2 in \cite{Kash}).
This curve has only one singular point which is a cusp with $N+1$ Newton pairs.
The decorations $\{(p_k,a_k)\}_{k=1}^{N+1}$
of the corresponding Eisenbud-Neumann splice diagram obtained computing the graph determinants are

\begin{equation*}\left\{ \begin{array}{l}
p_k= {\tilde n}_k \ \mbox{for $1\leq k\leq N,$} \ \mbox{and} \ p_{N+1}=\vpt^2;\\
a_k=(\vpo^2 {\tilde n}_1^2\cdots {\tilde n}_{k-1}^2 {\tilde n}_k-1)/{\vpt^2}, \ \mbox{for $1\leq k
\leq N$}\ \mbox{and} \ a_{N+1}=\vpo^2  {\tilde n}_1^2\cdots  {\tilde n}_{N}^2.
\end{array}\right.\tag{15} \end{equation*}

In the same way as before one checks that the $a_k$'s are integers.
The minimal set of generators of its semigroup $\Gamma_{II^-(\ell,N;\blambda)_{ge}}$ are
$\bar{\beta}_0=\vpt^2 {\tilde n}_1\cdots  {\tilde n}_N$,
$\bar{\beta}_k=(\vpo^2 n_1^2\cdots  {\tilde n}_{k-1}^2 {\tilde n}_k-1) {\tilde n}_{k+1}\cdots  {\tilde n}_N$ for $1\leq k\leq N$,
and $\bar{\beta}_{N+1}=\vpo^2  {\tilde n}_1^2\cdots  {\tilde n}_{N}^2.$

\subsubsection{$II^-(\ell,N;\lambda_1,\ldots,\lambda_N)_{sp}$}

The degree of the unicuspidal rational curve $II^-(\ell,N;\lambda_1,\ldots,\lambda_N)_{sp}$ is
$d=\vpo n_1\cdots n_N,$ (see Proposition 7.2 in \cite{Kash}).
Its singularity has $N+1$ Newton pairs and the decorations $\{(p_k,a_k)\}_{k=1}^{N+1}$
of the corresponding Eisenbud-Neumann splice diagram are

\begin{equation*}\left\{ \begin{array}{l}
p_k={\tilde n}_k \ \mbox{for $1\leq k \leq N,$} \mbox{and} \ p_{N+1}=\vpt;\\
a_k=(\vpo^2 {\tilde n}_1^2\cdots {\tilde n}_{k-1}^2{\tilde n}_k-1)/{\vpt^2}, \ \mbox{for $1\leq k
\leq N$}\ \mbox{and} \ a_{N+1}=(\vpo^2 {\tilde n}_1^2\cdots {\tilde n}_{N}^2+1)/{\vpt}.
\end{array}\right.\tag{16} \end{equation*}

Again $a_{N+1}$ is integer and the generators of its
semigroup $\Gamma_{II^-(\ell,N;\blambda)_{sp}}$ are
$\bar{\beta}_0=\vpt {\tilde n}_1\cdots {\tilde n}_N$,
$\bar{\beta}_k=(\vpo^2 {\tilde n}_1^2\cdots {\tilde n}_{k-1}^2{\tilde n}_k-1)
{\tilde n}_{k+1}\cdots {\tilde n}_N/{\vpt}$ for $1\leq k\leq N$, and
$\bar{\beta}_{N+1}=(\vpo^2 {\tilde n}_1^2\cdots {\tilde n}_{N}^2+1)/{\vpt}.$

\subsection{$(CP)$ for generic members of the pencils.}
The generic members of the pencils are the curves
$II^{\varepsilon}(\ell,N;\lambda_1,\ldots,\lambda_N)_{ge}, \, \varepsilon=\pm$.
To deal with the elements of their semigroups $\Gamma_{II^\pm(\ell,N;\blambda)_{ge}}$
it is better to multiply them by $\vpt.$
Define
\begin{equation*}
\chi^\varepsilon_{ge}(t):=\sum _{k\in\Gamma_{II^\varepsilon(\ell,N;\blambda)_{ge}}}\, t^{\lceil \vpt k/d\rceil} \ \mbox{and} \, \psi^\varepsilon_{ge}(t):=\frac{\chi^\varepsilon_{ge}(t) (1-t^{\vpt})}{1-t},  \ \mbox{with $\varepsilon=\pm$,}
\end{equation*}
where either $d=\vpt\vpf n_1\cdots n_N$ if $\varepsilon=+$ or $d=\vpt\vpo {\tilde n}_1\cdots {\tilde n}_N$ otherwise.
If
\begin{equation*}
CP^\varepsilon_{ge}(t):=\sum _{k\in\Gamma_{II^\varepsilon(\ell,N;\blambda)_{ge}}}\, t^{\lceil  k/d\rceil} \ \mbox{with $\varepsilon=\pm$,}
\end{equation*}
then \ref{compare} implies
\begin{equation*}
CP^\varepsilon_{ge}(t^{\vpt})=\frac{1}{\vpt} \sum_{\xi^{\vpt}=1} \psi^\varepsilon_{ge}(\xi t) \ \mbox{with $\varepsilon=\pm$.} \tag{17}
\end{equation*}

\subsection{Proposition.}\label{6.4} {\em

a) For the curve $II^{+}(\ell,N;\lambda_1,\ldots,\lambda_N)_{ge} $
the following identity holds:
\begin{equation*}
\chi^+_{ge}(t)=\frac{(1-t^{\vpo\vpf})(1-t^{\vpt^2\vpf n_1\cdots n_N})}{(1-t^{\vpt^2})(1-t^{\vpo})(1-t^{\vpf})}.
\end{equation*}

b) For the curve $II^{-}(\ell,N;\lambda_1,\ldots,\lambda_N)_{ge} $ the following identity holds:
\begin{equation*}
\chi^-_{ge}(t)=\frac{(1-t^{\vpo\vpf})(1-t^{\vpt^2\vpo {\tilde n}_1\cdots {\tilde n}_N})}{(1-t^{\vpt^2})(1-t^{\vpo})(1-t^{\vpf})}.
\end{equation*}}

\noindent {\em Proof.} We start with the curve
$II^{+}(\ell,N;\lambda_1,\ldots,\lambda_N)_{ge}$ which
has degree $d=\vpt\vpf n_1\cdots n_N$ and whose generators of its semigroup
$\Gamma_{II^+(\ell,N;\blambda)_{ge}}$ has been described in (12).
By \ref{a1}, each element $x$ in the semigroup can be written in a unique way as
$$
x=x_0 \vpt^2 \prod_{i=1}^N n_i+\sum_{k=1}^N y_k\ a_k\ \vpt^2 n_{k+1} \cdots n_N+ z_0\ \vpf^2\ n_1^2\cdots n_N^2,
$$
with $x_0\geq 0$, $0\leq y_k\leq n_k-1,$ for $1\leq k\leq N,$ and $0\leq z_0\leq \vpt^2-1.$
Thus
$$
\frac{\vpt x}{d}=x_0\frac{\vpt^2}{\vpf}+\sum_{k=1}^N y_k \frac{(\vpf^2 n_1^2 \cdots n_{k-1}^2n_k-1)}{\vpf n_1 \cdots n_k}+ z_0 \vpf n_1\cdots n_N.
$$
Write $x_0=q_0 \vpf+r_0 $ with $0\leq r_0\leq \vpf-1.$ Using $\vpt^2=\vpf\vpo-1$, we get
$$
\lceil \vpt x/d\rceil =q_0\vpt^2+r_0\vpo+ \sum_{k=1}^N y_k \vpf n_1 \cdots n_{k-1}+ z_0 \vpf n_1\cdots n_N,
$$
because
$
-1< (-r_0 n_1\cdots n_N-\sum_{k=1}^N y_kn_{k+1}\cdots n_N)/{\vpf n_1 \cdots n_N}\leq 0,
$
for every $0\leq y_k\leq n_k-1$, for $1\leq k\leq N,$ and $0\leq r_0\leq \vpf-1.$
The result is proved  because $
\chi^+_{ge}(t)$ is equal to
$$
\sum_{q_0\geq 0} \sum_{r_0=0}^{\vpf-1} \sum_{k=1}^N \sum_{y_k=0}^{n_k-1}\sum_{z_0=0}^{\vpt^2-1} t^{*}=\frac{(1-t^{\vpo\vpf})\prod_{k=1}^N(1-t^{\vpf n_1\cdots n_k})(1-t^{\vpt^2\vpf n_1\cdots n_N})}{(1-t^{\vpt^2})(1-t^{\vpo}) \prod_{k=1}^N (1-t^{\vpf n_1\cdots n_{k-1}}) (1-t^{\vpf n_1\cdots n_N})},
$$
where $*=q_0\vpt^2+r_0\vpo+ \sum_{k=1}^N y_k \vpf n_1 \cdots n_{k-1}+ z_0 \vpf n_1\cdots n_N.$

The proof for the curve $II^{-}(\ell,N;\lambda_1,\ldots,\lambda_N)_{ge}$ is essentially the same replacing
$\vpf$ by $\vpo.$

\vspace{2mm}

\subsection{Corollary.}{\em $(CP)$ is true for $II^{+}(\ell,N;\lambda_1,\ldots,\lambda_N)_{ge}$ and $II^{-}(\ell,N;\lambda_1,\ldots,\lambda_N)_{ge}$.}

\vspace{2mm}

\noindent {\em Proof.}
We do the case $\varepsilon=+$ leaving the other one to the reader.
Using \ref{6.4} and (17) then
\begin{equation*}
CP^+_{ge}(t^{\vpt})=\frac{(1-t^{\vpt})(1-t^{\vpt^2\vpf n_1\cdots n_N})}{(1-t^{\vpt^2})}\cdot \frac{1}{\vpt} \sum_{\xi^{\vpt}=1} \frac{(1-(\xi t)^{\vpo\vpf})}{(1-(\xi t)^{\vpo})(1-(\xi t)^{\vpf})(1-\xi t)}.
\end{equation*}
Let $\Delta(t)$ be the characteristic polynomial of the rational unicuspidal curve of
Example (d) in \ref{4.2}. That is, it has degree $\vpt$
and one characteristic pair $(a,b)=(\vpo,\vpf)$.
For such a curve we have already proved $R(t)\equiv 0$
in \ref{4.4}. Thus
$$
\frac{1}{\vpt}\sum_{\xi^{\vpt}=1} \frac{(1-(\xi t)^{\vpo\vpf})}{(1-(\xi t)^{\vpo})(1-(\xi t)^{\vpf})(1-\xi t)}=\frac{1}{\vpt}\sum_{\xi^{\vpt}=1} \frac{\Delta(\xi t)}{(1-\xi t)^2}=\frac{(1-t^{\vpt^2})}{(1-t^{\vpt^2})^3}.$$
Since the degree of $II^{+}(\ell,N;\lambda_1,\ldots,\lambda_N)_{ge}$ is $\vpt\vpf n_1\cdots n_N$ then
the result follows because
$$
CP^+_{ge}(t^{\vpt})=\frac{(1-t^{\vpt^2\vpf n_1\cdots n_N})}{(1-t^{\vpt^2})^2}.
$$

\vspace{2mm}

\subsection{$(CP)$ for special members of the pencils.}
The
special members of the pencils are the rational curves
$II^{\varepsilon}(\ell,N;\lambda_1,\ldots,\lambda_N)_{sp}, \, \varepsilon=\pm$.
As in the previous case, to deal with the elements of their semigroups $\Gamma_{II^\pm(\ell,N;\blambda)_{sp}}$
we multiply them by $\vpt.$
Define
\begin{equation*}
\chi^\varepsilon_{sp}(t):=\sum _{k\in\Gamma_{II^\varepsilon(\ell,N;\blambda)_{sp}}}\, t^{\lceil \vpt k/d\rceil} \ \mbox{and} \, \psi^\varepsilon_{sp}(t):=\frac{\chi^\varepsilon_{sp}(t) (1-t^{\vpt})}{1-t},  \ \mbox{with $\varepsilon=\pm$,}
\end{equation*}
where either $d=\vpf n_1\cdots n_N$ if $\varepsilon=+$ or $d=\vpo {\tilde n}_1\cdots {\tilde n}_N$ otherwise.
If
\begin{equation*}
CP^\varepsilon_{sp}(t):=\sum _{k\in\Gamma_{II^\varepsilon(\ell,N;\blambda)_{sp}}}\, t^{\lceil  k/d\rceil} \ \mbox{with $\varepsilon=\pm$,}
\end{equation*}
then \ref{compare} implies
\begin{equation*}
CP^\varepsilon_{sp}(t^{\vpt})=\frac{1}{\vpt} \sum_{\xi^{\vpt}=1} \psi^\varepsilon_{sp}(\xi t) \ \mbox{with $\varepsilon=\pm$.} \tag{18}
\end{equation*}

We do the proof for the curve $II^{+}(\ell,N;\blambda)_{sp}$ and the proof for the curve $II^{-}(\ell,N;\lambda_1,\ldots,\lambda_N)_{sp}$ is essentially the same replacing $n_i$ by ${\tilde n}_i$ and
$\vpf$ by $\vpo.$

The degree of $II^{+}(\ell,N;\lambda_1,\ldots,\lambda_N)_{g}$ is $d=\vpf n_1\cdots n_N$ and the generators of the semigroup $\Gamma_{II^+(\ell,N;\blambda)_{sp}}$
have been described in (13).
Each element $x$ in $\Gamma_{II^+(\ell,N;\blambda)_{sp}}$ can be written in a unique way as
$$
x=x_0 \vpt \prod_{i=1}^N n_i+\sum_{k=1}^N y_k\ a_k\ \vpt\ n_{k+1} \cdots n_N+ z_0\frac{\vpf^2 n_1^2\cdots n_N^2+1}{\vpt},
$$
with $x_0\geq 0,$ $0\leq y_k\leq n_k-1$, for $1\leq k\leq N,$ and $0\leq z_0\leq \vpt-1$.
Thus
$$
\frac{\vpt x}{d}=x_0\frac{\vpt^2}{\vpf}+\sum_{k=1}^N y_k \frac{(\vpf^2 n_1^2 \cdots n_{k-1}^2n_k-1)}{\vpf n_1 \cdots n_k}+ z_0 \frac{\vpf^2 n_1^2\cdots n_N^2+1}{\vpf n_1\cdots n_N}.
$$
Write $x_0=q_0 \vpf+r_0 $ with $0\leq r_0\leq \vpf-1$. Since $\vpt^2=\vpf\vpo-1$ then
$$
\lceil \vpt x/d\rceil =q_0\vpt^2+\alpha(r_0,{\mathbf y},z_0)+\lceil
\beta(r_0,{\mathbf y},z_0)\rceil.
$$
where $
\beta(r_0,{\mathbf y},z_0):=(z_0-r_0 n_1 \cdots n_N -\sum_{k=1}^N y_k  n_{k+1} \cdots n_{N})/{\vpf n_1\cdots n_N}
$ and
$$
\alpha(r_0,{\mathbf y},z_0):=r_0\vpo+ \sum_{k=1}^N y_k \vpf n_1 \cdots n_{k-1}+ z_0 \vpf n_1\cdots n_N.
$$
Set
$
A:=\{(r_0,{\mathbf y},z_0): 0\leq r_0\leq \vpf-1,\,0\leq z_0\leq \vpt-1,\,0\leq y_k\leq n_k-1,\, \mbox{for} \,k=1,\ldots,N\}.
$
Since for $(r_0,y_1,\ldots,y_N,z_0)\in A,$ one gets $-1<\beta(r_0,{\mathbf y},z_0)<1$ then $\lceil\,\beta(r_0,{\mathbf y},z_0) \rceil$ is either $0$ or $1.$
Consider the subset $P$ of $A$ defined by  $\{ z_0-r_0 n_1 \cdots n_N -\sum_{k=1}^N y_k  n_{k+1} \cdots n_{N}>0\}$
and $Q$ its complement in $A$. Thus $\lceil\,\beta(r_0,{\mathbf y},z_0) \rceil$ is $1$ if and only if $(r_0,{\mathbf y},z_0)\in P$
and it is zero otherwise. This implies
\begin{equation*}
\chi^+_{sp}(t):=\frac{\sum_{(r_0,{\mathbf y},z_0)\in P} t^{\alpha(r_0,{\mathbf y},z_0)+1}+\sum_{(r_0,{\mathbf y},z_0)\in Q} t^{\alpha(r_0,{\mathbf y},z_0)}}{(1-t^{\vpt^2})}.
\end{equation*}

\subsubsection{Remark.}\label{remark1}  If $(r_0,{\mathbf y},z_0)\in P$ then $0< z_0-r_0 n_1 \cdots n_N -\sum_{k=1}^N y_k  n_{k+1} \cdots n_{N}<\vpt.$

\noindent By (18), the following identity holds
$$CP^+_{sp}(t^{\vpt})=\frac{(1-t^{\vpt})}{(1-t^{\vpt^2})}\cdot\frac{1}{\vpt} \sum_{\xi^{\vpt}=1}\frac{\sum_{(r_0,{\mathbf y},z_0)\in P} (\xi t)^{\alpha(r_0,{\mathbf y},z_0)+1}+\sum_{(r_0,{\mathbf y},z_0)\in Q} (\xi t)^{\alpha(r_0,{\mathbf y},z_0)}}{(1-\xi t)}.$$

\noindent {\bf Claim:} the following identity holds too:
$$CP^+_{sp}(t^{\vpt})=\frac{(1-t^{\vpt})}{(1-t^{\vpt^2})}\cdot\frac{1}{\vpt} \sum_{\xi^{\vpt}=1}\frac{\sum_{(r_0,{\mathbf y},z_0)\in P} (\xi t)^{\alpha(r_0,{\mathbf y},z_0)}+\sum_{(r_0,{\mathbf y},z_0)\in Q} (\xi t)^{\alpha(r_0,{\mathbf y},z_0)}}{(1-\xi t)}.$$

\noindent For each $(r_0,{\mathbf y},z_0)\in P$ consider the series $t^{\alpha(r_0,{\mathbf y},z_0)+1}/(1-t)=t^{\alpha(r_0,{\mathbf y},z_0)+1}(1+t+t^2+..).$ Then
$$\frac{1}{\vpt} \sum_{\xi^{\vpt}=1}\frac{(\xi t)^{\alpha(r_0,{\mathbf y},z_0)+1}}{(1-\xi t)}$$
is nothing but the coefficients whose power is a multiple of $\vpt$ in $t^{\alpha(r_0,{\mathbf y},z_0)+1}(1+t+t^2+..).$ Therefore
to prove the {\bf claim} it is enough to check that for every $(r_0,{\mathbf y},z_0)\in P$ then $\alpha(r_0,{\mathbf y},z_0) \ne 0\,\mod \vpt.$ Write $\alpha(r_0,{\mathbf y},z_0)$ as
$$
r_0(\vpt-\vptw)+\left(\sum_{k=1}^N y_k n_1 \cdots n_{k-1}+ z_0 n_1\cdots n_N \right)(\vptw+2\vpt).
$$
Since $\gcd(\vptw,\vpt)=1$ (see property (iii) of the Fibonacci numbers) then $\alpha(r_0,{\mathbf y},z_0) = 0\,\mod \vpt$ if and only if
$$
-r_0+\sum_{k=1}^N y_k n_1 \cdots n_{k-1}+ z_0 n_1\cdots n_N = 0\,\, \mod \vpt.
$$
Moreover, since for all $k,$ $n_k =-1\,\, \mod \vpt$ then
$$
-r_0+\sum_{k=1}^N y_k (-1)^{k-1}+ z_0 (-1)^N = 0\,\, \mod \vpt\, \iff -(-1)^Nr_0-\sum_{k=1}^N y_k (-1)^{N-k}+ z_0 = 0\,\, \mod \vpt.
$$
But this last identity implies $z_0-r_0 n_1 \cdots n_N -\sum_{k=1}^N y_k  n_{k+1} \cdots n_{N}=s \vpt$ for $s\in \bz$ which contradicts the definition of the set $P$, see also \ref{remark1}.

\vspace{2mm}

\noindent Since
\begin{equation*}
\frac{\sum_{(r_0,{\mathbf y},z_0)\in P} t^{\alpha(r_0,{\mathbf y},z_0)}+\sum_{(r_0,{\mathbf y},z_0)\in Q} t^{\alpha(r_0,{\mathbf y},z_0)}}{(1-t^{\vpt^2})}=\frac{(1-t^{\vpo\vpf})(1-t^{\vpt\vpf n_1\cdots n_N})}{(1-t^{\vpt^2})(1-t^{\vpo})  (1-t^{\vpf})},
\end{equation*}
the {\bf claim} implies that $CP^+_{sp}(t^{\vpt})$ can be written as:
\begin{equation*}
\frac{(1-t^{\vpt})(1-t^{\vpt\vpf n_1\cdots n_N})}{(1-t^{\vpt^2})}\cdot \frac{1}{\vpt} \sum_{\xi^{\vpt}=1} \frac{\Delta(\xi t)}{(1-\xi t)^2}=
\frac{(1-t^{\vpt\vpf n_1\cdots n_N})}{(1-t^{\vpt})^2}
\end{equation*}
where again $\Delta(t)$ is the characteristic polynomial of the rational unicuspidal curve of degree $\vpt$
which has one characteristic pair $(a,b)=(\vpo,\vpf)$, see (d) in \ref{4.2}.
Since the degree of $II^{+}(\ell,N;\lambda_1,\ldots,\lambda_N)_{sp}$ is $\vpf n_1\cdots n_N$
then the result is proved.

\section{The case $\bar{\kappa}(\bp^2\setminus C)=1$ and
$\nu=1$.}\label{k=1n=1}

\subsection{Tono's classification theorem.}\label{Tsunoda}
We recall the following
classification result of K. Tono  \cite{KeitaTono}. Let $[x:y:z]$ be   a
system of homogeneous coordinates in $\bp^2$.
Assume that $C$ is a unicuspidal rational plane curve  with
$\bar{\kappa}=1$. Then $C$ is projectively equivalent with one of the
following curves $C'$:

\vspace{2mm}

\noindent {\em Type I.} $C'$ is given by
$$((f_1^sy+\sum_{i=2}^{s+1}a_if_1^{s+1-i}x^{ia-a+1})^{a}-f_1^{as+1})/x^{a-1}=0,$$
where $f_1=x^{a-1}z+y^a$, $a\geq 3$, $s\geq 1$, $a_2,\cdots, a_{s+1}\in \bc$
with  $a_{s+1}\not=0$.

In this case, $d=a^2s+1$, and the multiplicity sequence of $(C',p)$
is $[(a^2-a)s,(sa)_{2a-1},a_{2s}]$.

\vspace{1mm}

\noindent {\em Type IIa.} $C'$ is given by
$$((yf_2^n+x^{2n+1})^{4n+1}-f_3^{2n+1})/f_2^n=0,$$
where $f_2=xz-y^2$, $f_3=f_2^{2n}z+2x^{2n}yf_2^n+x^{4n+1}$ and
$n\geq 2$.

In this case $d=8n^2+4n+1$ and the multiplicity sequence of
$(C',p)$ is $[(n(4n+1))_4, (4n+1)_{2n},3n+1,n_3]$.

\vspace{1mm}

\noindent {\em Type IIb.} $C'$ is given by
$$((f_3^{2s-1}(f_2^ny+x^{2n+1})+\sum_{i=1}^sa_if_3^{2(s-i)}f_2^{i(4n+1)-n})^{4n+1}-f_3^{2((4n+1)s-n)})/f_2^n=0,$$
where $n\geq 2$, $s\geq 1$, $a_1,\ldots, a_s\in\bc$ with
$a_s\not=0$.

The degree of $C'$ is  $d=2(4n+1)^2s-4n(2n+1)$. Set $a^*:=4n+1$
and $s^*:=4s-1$. The multiplicity sequence for $s=1$ is
$[(3na^*)_4, (3a^*)_{2n}, (a^*)_3, 3n+1, n_3]$, otherwise it is
$$[(s^*a^*n)_4, (s^*a^*)_{2n}, (sa^*)_3, (s-1)a^*,
(a^*)_{2(s-1)}, 3n+1, n_3].$$

\subsection{$(C,p)$ in the case Type I}\label{top}
After a computation, one has the following facts:

If $s=1$, then $(C',p)$ has two characteristic pairs.
The decorations of the splice diagram
(for notation see \ref{ap}) are the following:
$$ p_1=a-1,\ p_2=a,\ a_1=a, \ \mbox{and} \
a_2=a(d+1)+1=a^3+2a+1.$$
In particular, the semigroup $\Gamma$ is generated by the elements
$a^2-a,\ a^2 $ and $a^3+2a+1$.

If $s>1$ the $(C',p)$ has three pairs: $p_1=a-1, \ p_2=s,\
p_3=a,\ a_1=a,\ a_2=d$ and $a_3=as(d+1)+1$. The semigroup is generated by
$$\bar{\beta}_0=(a^2-a)s,\ \bar{\beta}_1=a^2s,\ \bar{\beta}_2=
a(a^2s+1)=ad, \ \bar{\beta}_3=1+as(d+1)=a_3.$$ Notice that if in
this second case we consider $s=1$, we get a non-minimal splice
diagram (with $p_2=1$) of the first case. In particular, any
argument for general $s>1$ can be adopted to the $s=1$ situation
by a simple substitution $s=1$ and by elimination of the
semigroup generator $\bar{\beta}_2=ad$. Therefore, in the sequel
we will consider the general $s\geq 1$ situation  (which can be
specialized  to $s=1$ as described above).

\subsection{Theorem.}\label{00} {\em Type I satisfies the
distribution property ($CP$) (cf. \ref{fact3}): $R(t)\equiv0$.}

\begin{proof}
We will use the notations of \ref{fact3}. We will prove the inequality
\begin{equation*}
\#\Gamma \cap I_{l+1}\leq 1+\#\Gamma \cap I_{l},\tag{$in_l$}\end{equation*}
for any $l$. Notice that this implies  the negative distribution
property, hence the vanishing of $R(t)$ would  follow from \ref{nuone}.

Since $\Gamma\cap I_0=\{0\}$ and $\Gamma\cap I_1=\{(a^2-a)s,a^2s\}$,
the inequality ($in_0$) follows.
Since $d(d-3)+1$ is the largest element in $\bn\setminus \Gamma$,
($in_l$) for $l\geq d-2$ also follows. Next, we fix an $l$ with
$1\leq l\leq d-3$.
Denote by $i_0:=(l-1)d+1$ and  consider the map
$$\phi_l:\Gamma\cap I_l\setminus \{i_0\}\to \Gamma\cap I_{l+1}$$
defined by $\phi_l(x)=x+\bar{\beta}_1=x+a^2s$.
Clearly, $\phi_l$ is injective.

The following facts can be verified by elementary computations
(using the short hints):

(a) Using \ref{a1} (especially the formula for $\Delta$)
one gets that any $\gamma\in \Gamma$  can be written in a unique way in
the form
\begin{equation*}\gamma=k\bar{\beta}_3+m\bar{\beta}_0+t\bar{\beta}_2+n\bar
{\beta}_1,\tag{19}
\end{equation*}
where $0\leq k<a$, $0\leq m<a$, $0\leq t<s$  and $n\geq 0$. Here $k$ is the remainder
of $\gamma$ modulo $a$.

(b) If in (19)  $n>0$ then $\gamma=\gamma'+a^2s$
for some $\gamma'\in \Gamma$. In particular,
$\Gamma\cap I_{l+1}\setminus im(\phi_l)$ has two types of elements:

\vspace{2mm}

\noindent {\em Type A.} $\gamma=k\bar{\beta}_3+m\bar{\beta}_0+t\bar{\beta}_2$
with  $0\leq k<a$, $0\leq m<a $, $0\leq t<s$; or

\vspace{2mm}

\noindent {\em Type B.} $\gamma=(l+1)d$.

\vspace{2mm}

Here we do not exclude the situation when some semigroup element has
both types. Let $S_A$ (resp. $S_B$) be the set of elements of type A
(resp. B).

(c) $\#S_A\leq 2$.

Indeed, assume that
$\gamma=k\bar{\beta}_3+m\bar{\beta}_0+t\bar{\beta}_2$ and
$\gamma'=k'\bar{\beta}_3+m'\bar{\beta}_0+t'\bar{\beta}_2$
are both elements of $S_A$. If $k'<k$
then $(a-1)(a^2-a)s+(s-1)ad\geq m'\bar{\beta}_0+t'\bar{\beta}_2>ld-
k'\bar{\beta}_3\geq -d+(k-k')\bar{\beta}_3
\geq -d +\bar{\beta}_3$, a contradiction.
Hence $k=k'$.
Next, consider the difference $\gamma-\gamma'=(m-m')\bar{\beta}_0+(t-t')da$.
Since $|\gamma-\gamma'|<d$ and $|m-m'|\bar{\beta}_0<(a-1)d$
one gets $|t-t'|ad<ad$, hence $t=t'$.
Therefore, if $\gamma'>\gamma$ then $\gamma'=\gamma+(a^2-a)s$
(since $2(a^2-a)s\geq d$).

(d) One has the following identity for any $\gamma$ of type A:
\begin{equation*}
k\bar{\beta}_3+m\bar{\beta}_0+t\bar{\beta}_2=
(kas+m+at)d+(k-m)(1+sa),\tag{20}\end{equation*}
with $|k-m|\leq a-1$.

\vspace{2mm}

If $S_A=\emptyset $, or $S_A=\{\gamma\}$ and $(l+1)d\not\in \Gamma$,
or $S_A=\{(l+1)d\}$ then obviously $(in_l)$ follows.  Hence we only have to
analyze the following two cases:

\vspace{2mm}

\noindent {\bf Case 1.} $\#S_A=2$.

\vspace{2mm}

Write $S_A=\{\gamma,\gamma'\}$ where $\gamma'=\gamma+(a^2-a)s$
(see the proof of
(c)). Set $\gamma=ld+r$ for some $r> 0$. Since $ld+r+(a^2-a)s=\gamma'\leq
ld+d$, one gets $r\leq as+1$.

\vspace{2mm}

\noindent {\bf Case 1a.} First we verify that the case $0<r\leq as$ cannot
occur. Indeed, assume that $r\leq as$ and
write $\gamma$ in the form (20). Here there are two
possibilities: either $kas+m+at=l$ and $k-m>0$, or
$kas+m+at=l+1$ and $k-m<0$. The first possibility is eliminated by
$r=(k-m)(as+1)>as\geq r$. Hence $kas+m+ta=l+1$, and
$d-r=(m-k)(1+as)$. If $m-k\leq a-2$ then
$a^2s+1-as\leq d-r\leq (a-2)(1+sa)$,
a contradiction. Finally, if $m-k=a-1$ then
necessarily $k=0$ and $m=a-1$, hence $\gamma'=t\bar{\beta}_2+(a-1)
\bar{\beta}_1$ is not of type A,
again a contradiction.

\vspace{2mm}

\noindent {\bf Case 1b.} Therefore, $r=as+1$. In other words, $\gamma=ld+as+1$
and $\gamma'=(l+1)d$. In particular, $\gamma'$ is of also of type B.
This implies that  $\#\Gamma\cap I_{l+1}=2+\#im (\phi_l)$.
Next, write $\gamma$ as in (20). Then an elementary computation shows that
$kas+m+at=l+1$ is not possible. Therefore $kas+m+at=l$ and $k-m=1$.
In particular,
$$i_0=m\bar{\beta}_3+(m+a-1)\bar{\beta}_0+t\bar{\beta}_2+(a(s-1)+1)\bar{\beta}_1,$$
hence $i_0\in \Gamma\cap I_l$. This implies $
\#\Gamma\cap I_{l}=1+\#im (\phi_l)$, hence $(in_l)$  follows with equality.

\vspace{2mm}

\noindent {\bf Case 2.} $S_A=\{\gamma\}$, $(l+1)d\in \Gamma$ and $\gamma\not=
(l+1)d$.

\vspace{2mm}

Our goal is to show that $i_0\in \Gamma$.
Since $(l+1)d$ is not of type A, it follows that $(l+1)d=\tilde{\gamma}+a^2s$
for some $\tilde{\gamma}\in \Gamma$ -- in fact, for  $\tilde{\gamma}=ld+1$.
Hence we have to verify the following fact:
If $\tilde{\gamma}=ld+1\in \Gamma$, and there exists $\gamma\not=(l+1)d$
of type A in $\Gamma\cap I_{l+1}$, then $i_0\in \Gamma$.

Write $\gamma=k\bar{\beta}_3+m\bar{\beta}_0+t\bar{\beta}_2$ with
$0\leq k<a$, $0\leq m<a $, $0\leq t<s$. Similarly, set some
representation
$\tilde{\gamma}=k'\bar{\beta}_3+m'\bar{\beta}_0+t'\bar{\beta}_2
+n'\bar{\beta}_1$. Here we do not impose any inequality for the 
coefficients $k',t',m',n'$; hence $k'$ is uniquely determined by
$\tilde{\gamma}$, but $m'$ and $n'$ not (because of the relation
$a(a^2-a)=a^2(a-1)$). Notice also that $i_0=\tilde{\gamma}-d$.
Therefore, the identity
$$d=\bar{\beta}_3+2\bar{\beta}_0-(1+sa)\bar{\beta}_1$$
shows that if there is a representation of $\tilde{\gamma}$ with
$k'\geq 1$ and some $m'\geq 2$ then $i_0\in \Gamma$.

First we verify that (in any representation of $\tilde{\gamma}$ as above)
$k'\geq 1$. Indeed, assume that $k'=0$. Then (taking $\tilde{\gamma}=ld+1$
modulo $a$) one gets that
$l=-1+La$ for some $1\leq L<
as$. Since $ld+1\leq \gamma\leq (l+1)d-1$, one has
\begin{equation*}
La(a^2s+1)-a^2s\leq
k\bar{\beta}_3+m\bar{\beta}_0+t\bar{\beta}_2
 \leq La(a^2s+1)-1.\tag{21}
\end{equation*}
This implies that
$$\Big\lceil \frac{
k\bar{\beta}_3+m\bar{\beta}_0+t\bar{\beta}_2
}{a(a^2s+1)}\Big\rceil=L.$$ But one also has $$(t+ks)a(a^2s+1)<
k\bar{\beta}_3+m\bar{\beta}_0+t\bar{\beta}_2\leq
(t+ks+1)a(a^2s+1).$$
Hence $L=t+ks+1$. But for this $L$ the left inequality of (21)
fails.

Next, we show that there exists a representation with $m'\geq 2$.
E.g., if in some ``bad'' representation of $\tilde{\gamma}$ one
has $m'=0$, then taking
$k'\bar{\beta}_3+t'\bar{\beta}_2+n'\bar{\beta}_1$ modulo $d$ one
gets $n'\equiv -1+k'(as+1)$ modulo $d$.  Since $1\leq k'\leq
a-1$, this implies that $n'\geq as$, hence $n'a^2$ can be
rewritten as $(n'-a+1)a^2+a(a^2-a)$. Similar argument works for
$m'=1$ as well.
\end{proof}

\subsection{Addendum to the proof of \ref{00}}\label{01} In the verification
of the conjecture for  Type II curves we will use the results
valid for the Type I curves. We will need the following
additional facts.

 Assume that $(C',p)$ is of Type I as above, and
we keep the notations of \ref{top} and \ref{00}. \ref{00} shows
that ($*$)\  $\#\Gamma\cap [0,ld\,]=(l+1)(l+2)/2$ for any $0\leq
l<d$. We first assume that $d\equiv 0$, $a\equiv 1$ and $s\equiv
-1\mod 4$.

Recall that any $\gamma\in \Gamma$ can be written in a unique way
as
$\gamma=k\bar{\beta}_3+m\bar{\beta}_0+t\bar{\beta}_2+n\bar{\beta}_1$.
It is easy to see that $\gamma\equiv -n\mod 4$.  Set
$\#_{l,i}:=\#\{\gamma\in\Gamma,\ \gamma\leq ld,\ n\equiv i\
\mbox{mod}\ 4\}$ for $i=0,1,2,3$ and $0\leq l < d$. Then,  using
the function $x\mapsto x+\bar{\beta}_1=x+d-1$, one gets
$$\#_{l,1}=\#_{l-1,2},\
\#_{l,2}=\#_{l-1,3},\ \#_{l,3}=\#_{l-1,0},\
\#_{l,0}=\#_{l-1,1}+A$$ for some $A\in\bn$. Here one uses that
$4|ld$.  Taking the sum over $i$ and using ($*$), one gets that
$A=l+1$. Hence $\#_{l,i}$ can be computed inductively. In
particular, for any $l< d$,
$\#_{l,0}=(l+1)+\#_{l-4,0}=(l+1)+(l-3)+\#_{l-8,0}=\cdots. $ More
precisely, if $l=4k+r$ (with $0\leq r\leq 3$), then
$$\#_{l,0}=(k+1)(2k+r+1)\ \mbox{for any $l < d$}.$$
E.g., this can be rewritten into
\begin{equation*}\#_{2l,0}=(l+1)(l+2)/2, \ \mbox{whenever
$2l<d$}.\tag{$\#_{2l}$}\end{equation*} The very last identity is
true even if $s=1$, $a\equiv 1$ and $d\equiv 2$ $\mod 4$. Indeed,
in this case $\bar{\beta}_2$ is eliminated,
$\bar{\beta}_0\equiv\bar{\beta}_3\equiv \bar{\beta}_1-1\equiv 0$
$\mod 4$, hence $\gamma\equiv n$ $\mod 4$. Therefore,
$\#_{l,1}=\#_{l-1,0}$ and $\#_{l,3}=\#_{l-1,2}$ similarly as above,
 since $ld$ is even. Although $\#_{l,2}=\#_{l-1,1}$ may not be true in
general, but for $l$ even it is true. Writing $\#_{l,0}=\#_{l-1,3}+A$,
the formula
for $(\#_{2l})$ works again by a similar argument as above.

\subsection{$(C',p)$ in the case Type II}\label{type2} One can
verify the following facts:

In case IIa one has two characteristic pairs, the splice
decorations are: $p_1=n,\ p_2=4n+1,\ a_1=4n+1$ and
$a_2=(2n+1)d-n(4n+1)$.

In case IIb one has three characteristic pairs, and the splice
decorations can be uniformly described for any $s\geq 1$. They
are: $p_1=n,\ p_2=s^*,\ p_3=a^*,\ a_1=a^*,\ a_2=d/2$ and
$a_3=a^*s^*d/2+(s-1)a^*+3n+1$.

In fact, the IIa case also can be considered as a
`specialization': if in the formulas of IIb (including the
formula of $d$ as well) we substitute $s=1/2$ then we get the
non-minimal splice diagram (with $p_2=s^*=1$) of the IIa case.
Hence in the sequel we will handle the case II uniformly (using
the formulas of IIb, where $s$ is any positive integer or 1/2).

It is easy to verify that for Type II the degree $d$ satisfies
$2d=(a^*)^2s^*+1$, and the generators of $\Gamma$  are
$$\bar{\beta}_0=\frac{(a^*)^2-a^*}{4}\, s^*;\ \
\bar{\beta}_1=2d-1=(a^*)^2s^*;\ \ \bar{\beta}_2=a^*d/2;\ \
\bar{\beta}_3=\frac{a^*s^*}{4}(2d+1)+\frac{1}{4}.$$ If $s=1/2$,
or $s^*=1$, then one should eliminate  $\bar{\beta}_2$.

\subsection{Theorem.}\label{22} {\em The Conjecture is true for
Type II curves.}

\begin{proof}
One can also think about the integers $a^*$ and $s^*$ as
parameters of a Type I curve ($a^*$, $s^*$ corresponding to $a$
and $s$) with semigroup $\Gamma^*$. Their other  Type I invariants
are (cf. \ref{top}):
$$d^*=(a^*)^2s^*+1;\
\bar{\beta}_0^*=((a^*)^2-a^*)s^*;\
\bar{\beta}_1^*=d^*-1=(a^*)^2s^*;\  \bar{\beta}_2^*=a^*d^*;\
\bar{\beta}_3^*=a^*s^*(d^*+1)+1.$$ Notice that $$d^*=2d,\
\bar{\beta}_0^*=4\bar{\beta}_0,\ \bar{\beta}_1^*=\bar{\beta}_1,\
\bar{\beta}_2^*=4\bar{\beta}_2;\ \bar{\beta}_3^*=4\bar{\beta}_3.$$
We fix an integer $l<d$ and we wish to determine
$$\#\{\gamma\in\Gamma:\,
\gamma=k\bar{\beta}_3+m\bar{\beta}_0+t\bar{\beta}_2+n\bar{\beta}_1\leq
ld\}.$$ Multiplying by 4, one gets that this equals
$$\#\{\gamma^*\in\Gamma^*:\,
\gamma^*=k\bar{\beta}^*_3+m\bar{\beta}^*_0+t\bar{\beta}^*_2+4n\bar{\beta}^*_1\leq
2ld^*\}.$$ If $s\geq 1$, then $d^*\equiv 0;\ a^*\equiv 1; \
s^*\equiv -1 \mod  4$. If $s=1/2$, then $d^*\equiv 2$ and
$a^*\equiv 1$ $\mod 4$, and $s^*=1$. Hence this wanted number is
$\#_{2l,0}$ computed in \ref{01} (notice that $2l<2d=d^*$), and
equals $(l+1)(l+2)/2$. Hence
$$\#\{\gamma\in\Gamma:\,
\gamma=k\bar{\beta}_3+m\bar{\beta}_0+t\bar{\beta}_2+n\bar{\beta}_1\leq
ld\}=(l+1)(l+2)/2$$ for any $s$.   In particular, conjecture
follows for Type II as well.
\end{proof}

\section{Lin-Zaidenberg curves: $\bar{\kappa}=1$ and $\nu=2$.}

In this section we use an equivalent definition of \emph{Lin-Zaidenberg type} curves given
in terms of automorphisms of the affine plane, (see section 4 for another definition).
Consider a coordinate system $x,y$ in $\bc^2$.
Take a pair $(p,q)$ of relatively prime numbers with $p>q$. Let $C_{p,q}$
be the curve of $\bc^2$ defined by the equation $y^p+x^q=0$.
Consider $\bp^2$ together with a projective reference $(X,Y,Z)$. We embed $\bc^2$
into $\bp^2$ declaring that the image of the embedding is the open
subset $U_{Z}$ defined by $Z\neq 0$ and that $(x,y)=(X/Z,Y/Z)$.
We will denote by the same symbol the curve $C_{p,q}$ and its compactification to $\bp^2$.
Let $\phi:\bc^2\to\bc^2$ be an algebraic automorphism.
The embedding of $\bc^2$ into $\bp^2$ allows us to view any automorphism of
$\bc^2$ as a birational transformation of $\bp^2$.

\subsection{Definition.}\label{8.1} A curve $C$ is said to be of \emph{Lin-Zaidenberg type},
(\emph{LZ type} for short) if it is equal to $\phi^{-1}(C_{p,q})$ for a certain automorphism $\phi$
viewed as birational transformation of $\bp^2$. (Taking the inverse $\phi^{-1}$
instead of $\phi$ is just for notational convenience).

For a rational bicuspidal plane curve, K. Tono \cite{To} has proved that $C$ is of
\emph{LZ type} if and only if $\bar{\kappa}(\bp^2\setminus C)=1.$

\subsection{Theorem.}\label{8.2} {\em $(CP)$ is satisfied by
any LZ type curve; in other words $R(t)\equiv 0$. }

\vspace{2mm}

Let $C$ be a \emph{LZ type} curve. The curve $C$ has two
singularities, one at the origin of $\bc^2$ and the other at
infinity. Let $L_1(t)$ and $L_2(t)$ denote respectively the
Poincar\'e series of their corresponding semigroups, see \ref{a1}.
From the very definition (4) of $R(t)$, the vanishing $R(t)\equiv
0$ is equivalent to
\begin{equation*}
\frac{1}{d}\sum_{\xi^d=1}L_1(\xi t)L_2(\xi t)=\frac{1-t^{d^2}}{(1-t^d)^3}. \tag{22}
\end{equation*}

The singularity at the origin is isomorphic to $y^p+x^q=0$,
therefore
\begin{equation*}
L_1(t)=\frac{(1-t^{pq})}{(1-t^p)(1-t^q)}.\tag{23}
\end{equation*}

The combinatorics of the minimal embedded resolution of the singularity at
infinity of $C$ are closely
 related to the combinatorics of the minimal resolution of
the indeterminacy of $\phi$ as a birational transformation of $\bp^2$, and to the combinatorics of the minimal embedded resolution of the singularity that
$C_{p,q}$ has at infinity.

We will need the facts on algebraic automorphisms of $\bc^2$ from section 4; we refer to~\cite{Bo1},~\cite{Bo2} for details and more precise statements.
Let $\pi:X\to\bp^2$ the minimal resolution of the indeterminacy of $\phi$. Let
$L:=\bp^2\setminus\bc^2$ denote the line at infinity and consider $\calA$ the dual graph of the normal crossing divisor $\pi^*L$.
Recall that we have a total order of the vertices of the
graph (this order is recovered combinatorially by successively contracting $-1$ vertices and adjusting weights at each step). Denote by $\calA^*$ the graph
$\calA$ minus its first vertex (the one corresponding to the line at infinity). The automorphism $\phi$ admits a unique decomposition
as $\psi\comp H$ where $H$ is an affine transformation of $\bc^2$ and $\psi$ is the automorphism of $\bc^2$ which is composition of $\pi^{-1}$ and the
successive contraction of all the irreducible components of $\pi^*L$ with self intersection $-1$. Moreover all the components except the one corresponding to the
last vertex of $A$ are contracted. For our purposes it is sufficient to consider automorphisms for
which $H$ is the identity; we will assume this in the sequel.
The automorphism $\phi$ admits a
decomposition as $\phi=\phi_r\comp ...\comp\phi_1$, where $\phi_i$ is the result considering
first the blowing ups whose exceptional divisors belong to the $i$-th floor and then contracting all possible components with self-intersection $-1$.
The graph associated to the resolution of indeterminacy of $\phi_i$ is the elementary graph of length $n_i$.

The singularity at infinity of $C_{p,q}$ is at the point $(1:0:0)$ of $\bp^2$ and is defined by the equation $y^p+z^{p-q}=0$. Let $\sigma:Y\to\bp^2$ the
composition of blowing ups giving its minimal embedded resolution of this singularity. By the irreductibility of $C_{p,q}$ at $(1:0:0)$ the blowing ups whose
composition gives rise to $\sigma$ are totally ordered. Let $E=\sigma^{-1}(1:0:0)$ be the exceptional divisor of $\sigma$; denote
by $\calG$ its dual graph. Each vertex of $\calG$ corresponds to an irreducible component of $E$, and hence to one of the blowing-ups whose composition is
$\sigma$, and therefore they are totally ordered. It is known that $\calG$ is a linear graph whose first vertex is one of the extremes. Although
the ordering given to the vertices is not the same that the order given by the linearity of the graph we will denote by $v_1$ and $v_2$ the two extremal
vertices of the graph (the vertex $v_1$ is the first vertex in the given order, but $v_2$ need not be the second). We decorate the graph adding an arrow to
the vertex
corresponding with the irreducible component of $E$.

We must distinguish two cases:

\textbf{Case 1}: the first indetermination point of $\phi^{-1}$ is different to the point at infinity of $C_{p,q}$. Let $\psi:Z\to\bp^2$ be the blowing up
process giving the minimal embedded resolution of the singularity that $C$ has at infinity. It is easy to
check that $\psi$ is the sequence of blowing ups needed to resolve the indeterminacy of $\phi$, followed by
the sequence of blowing ups providing the embedded resolution of the singularity of $C_{p,q}$ at
infinity. The divisor $\psi^*C$ can be decomposed as
\[\psi^*C=\tilde{C}+\sum_{i=1}^km_kE_k,\]
where $\tilde{C}$ is the strict transform of $C$ and $E_1$,..., $E_k$ the irreducible components of the exceptional divisor ordered by appearance.

In this case it is easy to see that the decorated dual graph $\calB$ of the exceptional divisor is obtained by joining with an edge the last vertex of
$\calA^*$ with
the first vertex $v_1$ of $\calG$. We give a total order to the vertices of the resulting graph are totally ordered by considering first the vertices of
$\calA^*$ and after the vertices of $\calG$ with the previously given order. This order coincides with the natural order obtained defined identifying the
vertices with the divisors $E_1$,...,$E_k$.

With this information it is easy to compute the multiplicity
sequence of the singularity, and to deduce from it the degree of
$C$ and the coefficients $m_1$,...,$m_k$. The degree is
$d:=\deg(C)=pn_1\cdots n_r.$ For computing the Poincar\'e series
$L_2$ we need to know the coefficients $m_i$ of the univalent and
trivalent vertices of $\calB$. We observe that the univalent
vertices are the first vertex of each floor of $\calA^*$ (which
in the picture are the univalent central vertices of the floors),
the second vertex of the first floor, and the vertex $v_2$ of
$\calG$. We have exactly one trivalent vertex for each floor in
addition to the last vertex of $\calG$, where the arrow
decorating the graph is attached. The Poincar\'e series of the
semigroup of the singularity at infinity obtained is:
\begin{equation*}
L_2(t)=\frac{(1-t^{n_1^2\cdots
n_r^2p^2-pq})\prod_{i=1}^r(1-t^{(n_1^2\cdots n_i^2-1)n_i\cdots
n_rp})}{(1-t^{n_1^2\cdots n_r^2p-q})(1-t^{n_1\cdots
n_rp})\prod_{i=1}^r(1-t^{(n_1^2\cdots n_i^2-1)n_{i+1}\cdots
n_rp})}.\tag{24}
\end{equation*}

\textbf{Case 2}: the first indetermination point of $\phi^{-1}$ is equal to the point at infinity of $C_{p,q}$. We can further assume that the singularities
at infinity of $\phi_r^{-1}(C_{p,q})$ and $C$ do not have the same resolution graph: if this were the case we could work with $\phi_r^{-1}(C_{p,q})$ instead
of $C_{p,q}$ and with $\phi_{r-1}\comp ...\comp\phi_1\comp H$ instead of $\phi$, and we would be in Case 1. Let $k$ be the only integer such that
$p-kq<q<p-(k-1)q$; then the condition that $\phi_r^{-1}(C_{p,q})$ and $C$ do not have the same resolution graph is equivalent to $n_r\geq k+1$.

Let $\pi':X\to\bp^2$ be the minimal blowing-up process resolving the indeterminacy of $\phi^{-1}$.
Viewing $\phi$ as a composition of a sequence of blowing ups followed by a sequence of contractions it
is clear that the dual graph of $(\pi')^*L$ can be naturally identified with the graph $\calA$, but with
a different ordering of the vertices (in particular the first vertices are those of the $r$-th floor of
$\calA$). Noticing that $n_r\geq k+1$, a comparison between the blowing up sequences $\pi'$ and $\sigma$
shows that their first $k+1$ blowing up processes are the same in both sequences. Call $\sigma'$ the
composition of these blowing-ups and let
W:=$\{w'_1,...,w'_{k+1}\}$ be the ordered set of vertices of $\calG$ corresponding to the
exceptional divisors of these blowing-ups. We can decompose $\sigma$ as
$\sigma=\sigma_2\comp\sigma_1$, where $\sigma_1$ is the composition of blowing-ups whose associated
vertex does not belong to $W$ and $\sigma_2$ is the composition of the remaining blowing-ups.
Suppose that $w'_i$ is the first of
the vertices of $W$ such that its associated exceptional divisor meets
 the strict transform of $C_{p,q}$ by $\sigma_2$. Clearly $w'_i$ can be regarded also as a vertex of
$\calA$. Denote by $\pi_1$ the composition of those blowing ups of $\pi$ such that the vertex associated
to their exceptional divisor is smaller or equal than $w'_i$ with the ordering in $\calA$ induced by
$\pi$.

Let $\psi:Z\to\bp^2$ be the blowing up process giving the minimal embedded resolution of the
singularity that $C$ has at infinity. Then it is easy to check that $\psi=\pi_1\comp\sigma_1$. Due to this
decomposition
the decorated dual graph $\calB$ of the exceptional divisor of $\psi$ can be obtained as follows:

Let $\calA'$ be the graph obtained from $\calA^*$ by deleting the
last $k$ vertices. If $k>1$ let $\calG'$ be the graph obtained
from $\calG$ by deleting the first $k-1$ vertices; if $k=1$ let
$\calG'$ be the graph obtained by deletion of the first vertex of
$\calG$. The graph $\calG'$ is linear. If $k=1$ then its first
vertex is a extreme of it, and is denoted by $w$. If $k>0$ both
its first and second vertices are extremes of it; in this case we
alter the ordering in the vertices of $\calG'$ by interchanging
the order of the first two vertices, and denote by $w$ the first
vertex of $\calG'$ with the altered order. The graph $\calB$ is
is obtained by {\em identifying} (not joining with an edge) the
last vertex of $\calA'$ with the vertex $w$ of $\calG'$. We give
a total order to the vertices of the resulting graph by
considering first the vertices of $\calA'$ and after the vertices
of $\calG'$ with the altered order. This order coincides with the
natural order obtained defined identifying the vertices with the
irreducible components of the exceptional divisor of $\psi$. The
univalent vertices of $\calB$ are the first vertex of each floor
of $\calA$ (in the picture these are the univalent central
vertices of the floors), the second vertex of the first floor,
and the extreme of $\calG'$ different from $w$. We have exactly
one trivalent vertex for each floor in addition to the last
vertex of $\calG'$, where the arrow decorating the graph is
attached. The degree of $C$ is $d=n_1\cdots n_rq$ and the
Poincar\'e series of the singularity  at infinity is:
\begin{equation*}
L_2(t)=\frac{(1-t^{n_1^2\cdots
n_r^2q^2-pq})\prod_{i=1}^r(1-t^{(n_1^2\cdots n_i^2-1)n_i\cdots
n_rq})}{(1-t^{n_1^2\cdots n_r^2q-p})(1-t^{n_1\cdots
n_rq})\prod_{i=1}^r(1-t^{(n_1^2\cdots n_i^2-1)n_{i+1}\cdots
n_rq})} .\tag{25}
\end{equation*}

\subsubsection{Remark} \label{inter}
The formulas for the degree and $L_2(t)$ obtained in Case 2 are the result of interchanging the role of $p$ and $q$ in the formulas obtained in Case 1.

Due to $(22)$, $(23)$, $(24)$, $(25)$ and \ref{inter}, we need to
prove the following fact. For any two coprime positive integers
$p$ and $q$ (without imposing $p>q$), for $n_1,\ldots,n_r$
positive integers and $d:= n_1\cdots n_rp$ the following identity
holds:
\begin{equation*}
\frac{1}{d}\sum_{\xi^d=1}\frac{(1-(\xi t)^{pq})(1-(\xi
t)^{n_1^2\cdots n_r^2p^2-pq}) \prod_{i=1}^r(1-(\xi
t)^{(n_1^2\cdots n_i^2-1)n_i\cdots n_rp})}{(1-(\xi
t)^p)(1-t^q)(1-(\xi t)^{n_1^2\cdots n_r^2p-q})(1-(\xi
t)^{n_1\cdots n_rp})\prod_{i=1}^r(1-(\xi t)^{(n_1^2\cdots
n_i^2-1)n_{i+1}\cdots n_rp})}=\frac{1-t^{d^2}}{(1-t^d)^3}.
\end{equation*}

\noindent Multiplying in both sides by $d(1-t^d)^3$ and extracting common factors of numerator and denominator we see that the last equality is equivalent to:
\begin{equation*}
(1-t^d)\sum_{\xi^d=1}[(\sum_{j=0}^{p-1}(\xi
t)^{jq})(\sum_{j=0}^{p-1}(\xi t)^{j(n_1^2\cdots
n_r^2p-q)})(\sum_{j=0}^{n_1\cdots n_r-1}(\xi t)^{jp})\centerdot
\tag{26}
\end{equation*}
\[\centerdot\prod_{i=1}^r(\sum_{j=0}^{n_i-1}(\xi t)^{j(n_1^2\cdots
n_{i-1}^2n_i-1)n_{i+1}\cdots n_rp})]=d(1-t^{d^2}).\] The left hand
side of (26) is equal to
\begin{equation*}
(1-t^d)\sum_{\xi^d=1}[\sum_{i,j=0}^{p-1}\sum_{l=0}^{n_1\cdots
n_r-1}\sum_{\alpha=1}^r\sum_{k_\alpha=0}^{n_\alpha-1}(\xi
t)^{iq+j(n_1^2\cdots
n_r^2p-q)+lp+\sum_{\alpha=1}^rk_\alpha(n_1^2\cdots
n_{\alpha-1}^2n_\alpha-1)n_{\alpha+1}\cdots n_rp}].\tag{27}
\end{equation*}
The exponent of $(\xi t)$ in the last expression can be written as
$
A(i,j,l,{\mathbf k})+dB(i,j,l,{\mathbf k}),$
for
\[A(i,j,l,{\mathbf k})=(i-j)q+lp-\sum_{\alpha=1}^rk_\alpha
n_{\alpha+1}\cdots n_rp,\ \ B(i,j,l,{\mathbf k})=jn_1\cdots
n_r+\sum_{\alpha+1}^rk_\alpha n_1\cdots n_{\alpha-1}.\] But
$\sum_{\xi^m=1}\xi^m=0$ for any integer $m$ non-divisible by $d$.
In the present case $m=A(i,j,l,{\mathbf k})$.

%
%
%
%

\subsubsection{Lemma }\label{clave}
\emph{If $A(i,j,l,{\mathbf k})$ is divisible by $d$ then $i=j$
and $l=\sum_{\alpha=1}^rn_{\alpha+1}\cdots n_r$.}

As $d=n_1\cdots n_rp$ we deduce that $p|A(i,j,l,{\mathbf k})$. It
immediately follows that $i=j$ and that
\[n_1\cdots n_r|l-\sum_{\alpha=1}^rk_\alpha n_{\alpha+1}\cdots n_r.\]
We claim that then $0=l-\sum_{\alpha=1}^rk_\alpha
n_{\alpha+1}\cdots n_r$. The claim is proved by induction on $r$.
For $r=1$ it is obvious since $|l-k_1|<n_1$. For the induction
step we write
\[l-\sum_{\alpha=1}^rk_\alpha n_{\alpha+1}\cdots
n_r=l-k_r+n_r\sum_{\alpha=1}^{r-1}k_\alpha n_{\alpha+1}\cdots
n_{r-1};\] as $n_1\cdots n_r$ divides the last quantity we have
$n_r|l-k_r$. Therefore we can express $l$ as $l=k_r+n_rl'$. Then
\[n_1\cdots n_{r-1}|l'-\sum_{\alpha=1}^{r-1}k_\alpha n_{\alpha+1}\cdots n_{r-1},\]
and we conclude by induction.

\vspace{2mm}

Therefore, by \ref{clave} and  a computation, the left hand side
of  (26) becomes
\begin{equation*}
\label{LHS3} (1-t^d)\cdot
d\sum_{j=0}^{p-1}\sum_{\alpha=1}^r\sum_{k_\alpha=0}^{n_\alpha-1}t^{jn_1^2\cdots
n_r^2p+\sum_{\alpha=1}^rn_1^2\cdots n_{\alpha-1}^2n_\alpha\cdots
n_rp}.
\end{equation*}
A computation shows that the result equals $d(1-t^{d^2})$, as
desired.

\section{$(CP)$ for Orevkov's curves}

\subsection{}\label{5.1} In this section we will consider the families
$\{C_{4k}\}$ and $\{C^*_{4k}\}$ ($k\geq 2$)  of Orevkov \cite{Or}.
There is a special interest in these curves: Orevkov proved that they
satisfy $d>\alpha m$ (where $\alpha=(3+\sqrt{5})/2$ as above);
moreover he conjectured that these are the only rational cuspidal curves
(together with the curves $C_4$ and $C_4^*$ from \ref{4.2}(e)-(f))
satisfying this inequality. Also, for these curves one has
 $\bar{\kappa}(\bp^2\setminus C)=2$, in contrast with the previous section
4-8. In particular,
 \ref{5.2} shows that $(CP)$ is not a specialty of curves with
$\bar{\kappa}<2$ (and of some  finitely many sporadic curves of general type).

In this case the number of characteristic pairs of $(C,p)$ is again two,
which makes the structure of the semigroups rather  interesting.
The numerical invariants of the curves are the following (cf. \cite{Or}).
In both cases $j\equiv 0$ mod 4, and $j\geq 8$.

\subsubsection{} $C_{j}$ has degree $d=\vp_{j+2}$.  The
germ $(C_j,p)$  has two characteristic pairs with numerical invariants:
$$p_1=\vp_j/3,\ p_2=3;\ a_1=\vp_{j+4}/3, \ a_2=1+\vp_j\vp_{j+4}/3;$$
$$\bar{\beta}_0=\vp_j,\ \bar{\beta}_1=\vp_{j+4}, \
\bar{\beta}_2=a_2.$$

\subsubsection{} $C^*_j$ has degree $d^*=2\vp_{j+2}$.  The
germ $(C_j,p)$  has two characteristic pairs with numerical invariants:
$$p_1=\vp_j/3,\ p_2^*=6;\ a_1=\vp_{j+4}/3, \ a_2^*=1+2\vp_j\vp_{j+4}/3;$$
$$\bar{\beta}_0^*=2\vp_j,\ \bar{\beta}_1^*=2\vp_{j+4}, \
\bar{\beta}_2^*=a_2^*.$$

\subsection{Theorem.}\label{5.2} {\em $C_{4k}$ and $C^*_{4k}$ satisfy
$(CP)$.}

\vspace{2mm}

\noindent {\em Proof.} We start with the case $C_j$.

Step 1. Assume that $l<d$. Notice that gcd$(
\bar{\beta}_1,d)=1$, hence $k\bar{\beta}_1=ld$ implies $l=k=0$.
In all the other cases, by a similar argument as in the proof of
\ref{4.4}, $k\bar{\beta}_1<ld$ is equivalent with $k/l<1/\alpha$.
Hence $k\bar{\beta}_1\in I_l$ if and only if $\lceil k\alpha\rceil =l$.

Step 2. Assume $l<d$.
If $i\bar{\beta}_0>ld$ then $i/l>\vp_{j+2}/\vp_j>\alpha$.
Conversely, if $i/l>\alpha$ then consider $\vp_{j+2}/\vp_j>\vp_{j+4}/\vp_{j+2}
>\alpha$. Similarly as in the proof of \ref{4.4}, we get that $i/l\geq
\vp_{j+2}/\vp_j$. Hence, if $l\not\in\{\vp_j,\vp_j+1\}$ then $i\bar{\beta}_0
\in I_l$ if and only if $\lceil i/\alpha \rceil =l$.

Step 3. Notice that $x:=\vp_j\vp_{j+4}/3=(\vp_{j+2}^2-1)/3$
is not a multiple of $d$,
hence $x$ and $a_2=x+1$ are in the same interval $I_{l_0}$.
Let $\Gamma_0$ be the semigroup generated by $\bar{\beta}_0$
and $\bar{\beta}_1$. Hence for $l<l_0$ one has $\Gamma_0\cap I_l=
\Gamma\cap I_l$.

Moreover, $l_0-1<\vp_j$. Indeed, $(l_0-1)\vp_{j+2}<x=(\vp_{j+2}^2-1)/3<
\vp_{j+2}\vp_j$ (because $\vp_{j+2}<3\vp_j$).

Notice also that the very first semigroup element $i\vp_j+k\vp_{j+4}$
of $\Gamma_0$ which can be written in two different ways is $x=(\vp_{j+4}/3)
\vp_j=(\vp_j/3)\vp_{j+4}$.

Hence, for any $l\leq l_0-1$, by similar argument as in the proof of
\ref{4.4}(d), and using Step 1 and Step 2, we get that the distribution
patern is true.

\vspace{1mm}

The proof now bifurcates in two cases.

First assume that $d=\vp_{j+2}=3t-1$ for some $t$. Then $x=3t^2-2t$ and
$l_0=t$.

Step 4. One can verify that the following numbers are in increasing order:
$x-2\vp_j$, $(l_0-1)d$, $x-\vp_j$, $x$, $l_0d$, $x+\vp_j$, $x+3\vp_j$,
$(l_0+1)d$, $x+4\vp_j$, $x+6\vp_j$, $(l_0+2)d$, $x+\vp_{j+4}$, $x+7\vp_j$.
E.g., $x+\vp_j>l_0d$
reduces to $\vp_j>t$, or to $3\vp_j>\vp_{j+2}+1$ which is true.
The other  verifications are similar.

Step 5. Since $\vp_{j+4}+\vp_j=3\vp_{j+2}$, the map
$s_l:I_{l-3}\cap \Gamma_0\to I_l\cap \Gamma_0$, given by $y\mapsto
y+\vp_{j+4}+\vp_j$ is well-defined. Clearly, it is injective.
For $l=l_0$ the complement of its image has two elements, namely $x$ and
$x-\vp_j$ (cf. Steps 3 and  4).

For $l>l_0$  the map $s_l$ is surjective. This follows from Step 4
and the identities $i\vp_j=(i-\vp_{j+4}/3)\vp_j+(\vp_j/3)\vp_{j+4}$
(for $i>\vp_{j+4}/3$) and
$k\vp_{j+4}=(k-\vp_{j}/3)\vp_{j+4}+(\vp_{j+4}/3)\vp_{j}$
(for $k>\vp_{j}/3$).
Therefore
$$\sum_{k\in \Gamma_0}t^{\lceil k/d\rceil}=1+2t+\cdots + l_0t^{l_0-1}+
\Big(l_0+(l_0-1)t+l_0t^2\Big)\cdot(t^{l_0}+t^{l_0+3}+t^{l_0+6}+\cdots).$$

Step 6. Consider now the intervals $l_0\leq l\leq 2l_0-1$. Since (cf. Step 4)
$2x+d>2x+2\vp_j>2l_0d$, we get that $2a_2$ is not situating in these
intervals.By a mod 3 argument, for any such $l$,
$\Gamma\cap I_l$ is the disjoint union of $I_l\cap \Gamma_0$
and $I_l\cap (a_2+\Gamma_0)$. Moreover, in all these intervals
any element of $a_2+\Gamma_0$ has a unique representation in the form
$a_2+i\vp_j+k\vp_{j+4}$.

In particular, we have to  understand the
distribution of $a_2+\Gamma_0$.
By a similar computation as in Step 4, we get that $I_{l_0}$ contains
only $a_2$, $I_{l_0+1}$  contains $a_2+i\vp_j$ for $i=1,2,3$;
$I_{l_0+2}$ contains $a_2+i\vp_j$ for $i=4,5,6$. In particular, for intervals
$l=l_0,l_0+1,l_0+2$ the distribution of $\Gamma$ follows.

Step 7. Assume again $l_0\leq l\leq 2l_0-1$.
Consider the injective map $s_l':I_{l-3}\cap (a_2+\Gamma_0)
\to I_l\cap (a_2+\Gamma_0)$ given by $y\mapsto y+\vp_j+\vp_{j+4}$.
We claim that for $l_0-3\leq l\leq 2l_0-1$ the complement of the image of
$s_l'$
has exactly three elements. There are two types of elements in the
complement. The first type is
$a_2+k\vp_{j+4}=1+(k+\vp_j/3)\vp_{j+4}=1+k'\vp_{j+4}$.
Since $k'\vp_{j+4}$ is never a multiple of $d$ (in the relevant intervals)
 we get (via Step 1) that $1+k'\vp_{j+4}\in I_l$ if and only if
$k'\vp_{j+4}\in I_l$ if and only if $\lceil k'\alpha\rceil =l$.

The other type is $a_2+i\vp_j=1+(i+\vp_{j+4}/3)\vp_j=1+i'\vp_j$.
Notice that $i'\vp_j$ can be a multiple of $d$ in these intervals,
namely for $i'=d$. This fact combined with Step 2 we get that
$1+i'\vp_j\in I_l$ if and only if $\lceil i'/\alpha \rceil =l$.
Now, by the end of the proof of  \ref{4.4}(d) the above claim follows.

In particular, the distribution is true for any $l<2l_0$.

Step 8. Since $(CP_l)$ is true if and only if
$(CP_{d-2-l})$ is true, one gets $(CP_l)$ for all the remaining cases.

\vspace{2mm}

The proof in the other case $d=\vp_{j+2}=3t+1$ is similar.
The only differences are the following. First one has the following
increasing numbers:

$x-\vp_j$, $(l_0-1)d$, $x$, $x+\vp_j$, $l_0d$, $x+2\vp_j$, $x+4\vp_j$,
$(l_0+1)d$, $x+5\vp_j$, $x+6\vp_j$, $(l_0+2)d$,  $x+7\vp_j$.
Also $x+\vp_{j+4}<(l_0+2)d$.
In particular,
$$\sum_{k\in \Gamma_0}t^{\lceil k/d\rceil}=1+2t+\cdots +l_0t^{l_0-1}+
\Big((l_0-1)+(l_0-1)t+l_0t^2\Big)\cdot(t^{l_0}+t^{l_0+3}+t^{l_0+6}+\cdots).$$
Moreover, $a_2+\Gamma_0$ has 2, 3, resp. 3 elements in $I_l$ for
$l=l_0,l_0+1$, resp. $l_0+2$, namely
$a_2+i\vp_j$ for $i=0,1$ ($l=l_0$); $a_2+i\vp_j$ for $i=2,3,4$ ($l=l_0+1$);
and finally $a_2+i\vp_j$ for $i=5,6$ and $a_2+\vp_{j+4}$ in the last case.
Otherwise all the arguments are similar.

\vspace{2mm}

In the next paragraph we show that the case $C^*_j$ can be reduced to
the case $C_j$.

Assume first that $l<\vp_{j+2}$. In these intervals $I_l(C_j)$,
for any element of the semigroup of $C_j$, the coefficient $c_2$
of $a_2$ is $\leq 2$. Indeed, $3a_2=3+\vp_j\vp_{j+4}=2+\vp_{j+2}^2$.
Moreover, by a computation one can verify that there no elements of
type $2a_2+i\vp_j+k\vp_{j+4}$ (i.e. with $c_2=2$)
which equals some $l\vp_{j+2}+1$ (the first entry of some interval).
Using this, one gets a bijection
$\Gamma(C_j)\cap I_l(C_j)\to \Gamma(C^*_j)\cap I_l(C^*_j)$
given by $x\mapsto 2x-c_2$ (which sends $\bar{\beta}_i$ into
$\bar{\beta}^*_i$). Therefore, $(CP_l)$ is true for $C_j^*$ for any
$l\leq d^*/2-1$. By symmetry, $(CP_l)$ is true for any
$l\geq d^*-2-(d^*/2-1)=d^*/2-1$, hence for any $l$.

\end{document}